\documentclass{article}
\usepackage{amsmath,amssymb,mathrsfs,amsthm,verbatim}

\newtheorem{Theorem1}{Theorem}
\newtheorem{Theorem2}[Theorem1]{Theorem}

\newtheorem{Lemma1}{Lemma}[section]
\newtheorem{Lemma2}[Lemma1]{Lemma}
\newtheorem{Local}[Lemma1]{Lemma}
\newtheorem{Lemma3}[Lemma1]{Lemma}
\newtheorem{Lemma4}[Lemma1]{Lemma}
\newtheorem{Lemma5}[Lemma1]{Lemma}
\newtheorem{Lemma6}[Lemma1]{Lemma}

\newtheorem{Corollary}[Lemma1]{Corollary}
\newtheorem{Lemma9}[Lemma1]{Lemma}
\newtheorem{Lemma10}[Lemma1]{Lemma}
\newtheorem{Lemma11}[Lemma1]{Lemma}
\newtheorem{Lemma12}[Lemma1]{Lemma}
\newtheorem{Lemma13}[Lemma1]{Lemma}

\newtheorem{acknowledgement}{Acknowledgement}

\date{}

\begin{document}

\title{Global Existence for the Seiberg-Witten Flow}
 \author{ {Min-Chun Hong and Lorenz Schabrun} \footnote{
Department of Mathematics, The University of Queensland, Brisbane,
QLD 4072, Australia. Email: hong@maths.uq.edu.au,
lorenz@maths.uq.edu.au}}

\maketitle \textbf{Abstract.} We introduce the gradient flow of
the Seiberg-Witten functional on a compact, orientable Riemannian
4-manifold and show the global existence of a unique smooth
solution to the flow. The flow converges uniquely  in $C^\infty$
up to gauge to a critical point of the Seiberg-Witten functional.

\section{Introduction}
In his ground-breaking  work, Donaldson applied Yang-Mills theory
to construct a new invariant for 4-manifolds and proved that
there exist topological 4-manifolds which do not admit smooth
structures, and topological 4-manifolds that admit an infinite
number of distinct smooth structures (e.g. \cite{Donaldson}). A decade later,
Seiberg and Witten, again using considerations from gauge theory,
produced some surprisingly simple equations which have been used to produce simpler
proofs of many results from Donaldson theory, and also some new
results \cite{Witten.monopoles}. In particular, the new equations are first order and have gauge group
$U(1)$. Because of its ease of computation, Seiberg-Witten theory has effectively
succeeded Donaldson theory in many cases.

Computing the Seiberg-Witten invariant for a given manifold
involves finding nontrivial solutions to the Seiberg-Witten
equations (\ref{sweq}), called Seiberg-Witten monopoles.
Therefore, an important problem in Seiberg-Witten theory is the
formulation of necessary and/or sufficient conditions for the
existence of monopoles. In \cite{Taubes.Symplectic}, for instance,
Taubes proved that when a symplectic structure exists on $M$,
there exists a monopole for a particular canonical spin$^c$
structure. For an elementary introduction to spin geometry and the
Seiberg-Witten functional, see \cite{Jost.Riemannian}. For a
longer exposition of Seiberg-Witten theory, see \cite{Moore},
\cite{Morgan}, \cite{nicolaescu} or \cite{wildworld}.

Let $M$ be a compact oriented Riemannian 4-manifold with a
spin$^c$ structure $\mathfrak s$. Denote by $\mathcal S= W\otimes
\mathcal L$  the corresponding spinor bundle and by $\mathcal S^
\pm= W^{\pm}\otimes \mathcal L$ the half spinor bundles, and by
$\mathcal L^2$ the corresponding determinant line bundle. Recall
that the bundle $\mathcal S^+$ has fibre $\mathbb{C}^2$. Let $A$
be a unitary connection on $\mathcal L^2$. Note that we can write
$A=A_0 + a$, where $A_0$ is some fixed connection and $a \in
i\Lambda ^1 M$ with $i=\sqrt {-1}$. Denote by $F_A = dA \in
i\Lambda ^2 M$ the curvature of the line bundle connection $A$.
Let $\{e_j\}$ be an orthonormal basis of $\mathbb R^4$.  A
$\operatorname{Spin(4)^c}$-connection on the bundles $\mathcal S$ and $\mathcal
S^ \pm$ is locally defined by
\begin{equation} \label{connectiondef}
\nabla _A  = d + \frac{1}{2}(\omega +A),
\end{equation}
where  $\omega=\omega_{jk}e_je_k$ is induced by the
Levi-Civita connection matrix  $\omega_{jk}$ and $e_j e_k$ acts by
Clifford multiplication (see \cite{Jost.Riemannian}). We denote
the curvature of $\nabla _A$ by $\Omega _A$. The Dirac operator
$D_A :\Gamma (S) \to \Gamma (S)$ is given by
\[
D_A \varphi  = e_j \nabla _A^j\varphi ,
\]
where $\nabla _A^j$ denotes covariant differentiation along the
tangent vector $e_j$, and $e_j$ acts via Clifford multiplication.
We define the configuration space $\Gamma (\mathcal S^ +  )
\times \mathscr A$, where $\mathscr A$ is the space of unitary
connections on $\mathcal L^2$, and let $(\varphi ,A) \in \Gamma
(\mathcal S^ +  ) \times \mathscr A$.

The Seiberg-Witten equations are
\begin{equation} \label{sweq}
D^+_A \varphi  = 0,\;\;\;\;\;\;\;\;\;\;\;\;\;F_A^ +   = \frac{1}
{4}\left\langle {e_j e_k \varphi ,\varphi } \right\rangle e^j
\wedge e^k.
\end{equation}
Solutions with $\varphi=0$ are called reducible (or trivial)
solutions. Nontrivial solutions are called (Seiberg-Witten)
monopoles.
\medskip

The heat flow for the Yang-Mills equations, suggested by Atiyah and Bott, has played an important role in Yang-Mills theory. The first contribution was made by Donaldson [D2] in the case of a holomorphic vector bundle. He used the Yang-Mills heat flow to establish an important relationship
between Hermitian Yang-Mills connections and stable holomorphic
vector bundles. How to formulate a heat flow for the Seiberg-Witten equations and use it to establish a relationship
between Seiberg-Witten monopoles and spinor bundles is a challenging question.

\medskip
In order to answer this question, we introduce the gradient flow of
the Seiberg-Witten functional. The Seiberg-Witten functional
$\mathcal{SW}:\;\Gamma (\mathcal S^ +  ) \times \mathscr A \to
\mathbb{R}$ is given by
\[
\mathcal{SW}(\varphi ,A) = \int_M {\left| {D_A \varphi } \right|^2
}  + \left| {F_A^ +   - \frac{1} {4}\left\langle {e_j e_k \varphi
,\varphi } \right\rangle e^j  \wedge e^k } \right|^2\,dV.
\]
Using the Weitzenb\"ock formula (e.g. \cite {Jost.Riemannian} or
\cite {Moore})
\begin{equation} \label{Weitzenbock}
D_A^2 \varphi =  -\nabla^*_A\nabla_A \varphi + \frac{S} {4}
\varphi + \frac{1} {4}F_{A, jk}(e_je_k  \varphi),
\end{equation}
the Seiberg-Witten functional can be written in the following
form:
\begin{equation} \label{swfunc}
\mathcal{SW}(\varphi ,A) = \int_M {\left| {\nabla _A \varphi }
\right|^2  + \left| {F_A^ +  } \right|^2  + \frac{S} {4}\left|
\varphi  \right|^2  + \frac{1} {8}\left| \varphi  \right|^4 }\,dV,
\end{equation}
where $S$ is the scalar curvature of $M$. The Seiberg-Witten
functional is invariant under the action of a gauge group. The group of gauge
transformations is
\[
\mathscr G = \left\{ {g:M \to U(1)} \right\}.
\]
$\mathscr G$ acts on elements of the configuration space via
\[
g^*(\varphi ,A) = (g^{-1}\varphi ,A + 2g^{-1}dg).
\]
It is easily seen that (\ref{sweq}) and (\ref{swfunc}) are
invariant under the action of the gauge group.

Using the relation
\begin{equation} \label{selfdualrelation}
\left\| {F_A } \right\|_{L^2 }  = 2\left\| {F_A^ +  }
\right\|_{L^2 }  - 4\pi ^2 c_1 (\mathcal L)^2,
\end{equation}
where $c_1 (\mathcal L)$ is the first Chern class of $\mathcal{L}$
(see \cite{wildworld}), one can also write the functional in the
form
\begin{equation} \label{swfuncchern}
\mathcal{SW}(\varphi ,A) = \int_M {\left| {\nabla _A \varphi }
\right|^2  + \frac{1}{2}\left| {F_A} \right|^2  + \frac{S}
{4}\left| \varphi  \right|^2  + \frac{1} {8}\left| \varphi
\right|^4 }\,dV + \pi ^2 c_1 (\mathcal L)^2.
\end{equation}
Note that the term $\pi ^2 c_1 (\mathcal L)^2$ is constant along
the flow and does not affect the flow equations. Thus in this
paper, it can usually be neglected. The Euler-Lagrange equations
for the Seiberg-Witten functional are

\begin{equation} \label{eq1}
 - \nabla _A^* \nabla
_A \varphi  - \frac{1}{4}\left[ {S + \left| \varphi  \right|^2 }
\right]\varphi =0,
\end{equation}
\begin{equation} \label{eq2}
 - d^* F_A   - i\operatorname{Im} \left\langle {\nabla_A \varphi ,\varphi }
\right\rangle =0.
\end{equation}
The Euler-Lagrange equations for the Seiberg-Witten functional
were first  investigated by Jost, Peng and Wang in
\cite{Jost.Variational}. They proved a number of properties
including the Palais-Smale condition, compactness and the
smoothness of weak solutions to the system (\ref{eq1})-(\ref{eq2}). Note that equations (\ref{eq1})-(\ref{eq2}) always
admit the trivial solutions with $\varphi =0$, but among
the solutions to (\ref{eq1})-(\ref{eq2})  are also any nontrivial
solutions, including the Seiberg-Witten monopoles (solutions of
(\ref{sweq})).

Given the above functional setting, the natural evolution equation to
choose for finding critical points is the gradient flow.
Therefore, we define the Seiberg-Witten flow by
\begin{equation} \label{flow1}
\frac{{\partial \varphi }} {{\partial t}} =  - \nabla _A^* \nabla
_A \varphi  - \frac{1}{4}\left[ {S + \left| \varphi  \right|^2 }
\right]\varphi,
\end{equation}
\begin{equation} \label{flow2}
\frac{{\partial A}} {{\partial t}} =  - d^* F_A   -
i\operatorname{Im} \left\langle {\nabla_A \varphi ,\varphi }
\right\rangle
\end{equation}
with initial data
\[
(\varphi (0),A(0)) = (\varphi _0 ,A_0).
\]

Note that since the connection $\nabla_A$ respects the splitting $\mathcal{S}= \mathcal{S^+} \oplus  \mathcal{S^-}$, for initial data $\varphi_0 \in \Gamma (S^+)$, we have $\varphi(t) \in \Gamma (S^+)$ for each $t$. In this paper we establish that these flow equations
admit a smooth solution for all time, which converges to a
critical point of the functional (\ref{swfunc}).

\begin{Theorem1} \label{Main1}
For any given smooth  $(\varphi_0, A_0)$, the system (\ref{flow1})-(\ref{flow2}) admits a unique global smooth solution on $M
\times \left[ {0,\infty } \right)$ with initial data $(\varphi_0,
A_0)$.
\end{Theorem1}

We show the existence of a local solution to (\ref{flow1}) and
(\ref{flow2}) following an idea of Donaldson for the Yang-Mills
flow (e.g. \cite {Donaldson}) which considers a gauge equivalent flow. The
critical question for the global existence of the  Seiberg-Witten
flow turns out to be whether or not the energy concentrates, as
in the Yang-Mills and Yang-Mills-Higgs flows (see
\cite{StruweYMLocal} and \cite{Hong.YMHlocalexist}).
 While this question remains
unresolved for the Yang-Mills and Yang-Mills-Higgs flows in four
dimensions (see e.g. \cite{Hong.equivariant}), we fortunately show
that concentration does not occur in general for the
Seiberg-Witten flow at any time $T\leq \infty$.

\medskip

Concerning the limiting behaviour of the flow, we show the
following theorem.

 \begin{Theorem2} \label{Main2}
As $t\to\infty$, the solution $(\varphi (t),A(t))$ converges
  smoothly, up to gauge transformations,  to a unique limit
$(\varphi_{\infty}, A_{\infty})$, where $(\varphi_{\infty},
A_{\infty})$ is  a smooth solution of equations
(\ref{eq1})-(\ref{eq2}). There are  constants $C_k$ and $\frac 1 2
<\gamma <1$ such that
\begin{equation} \label{finalestimate}
\left\| {(\varphi (t),A(t)) - (\varphi _\infty  ,A_\infty )} \right\|_{H^k }  \leqslant C_kt^{ - (1 - \gamma )/(2\gamma  - 1)}.
\end{equation}
Moreover,  for any $\lambda
>0$, $(\varphi_0, A_0)\to (\varphi_{\infty}, A_{\infty})$ defines
a continuous map on the space  $\{(\varphi_0, A_0): SW  (\varphi
(t),A(t)) \to \lambda \}$ as $t\to \infty$.
\end{Theorem2}
Analogous results were proven for the Yang-Mills flow in two and
three dimensions by R\aa de \cite{Rade}, and extended to the
Yang-Mills-Higgs functional on a Riemann surface by Wilkin
\cite{Wilkin}. Both of these extend the work of Simon \cite{Simon}.

\medskip
Let $\mathcal A = \Gamma (\mathcal S^ +  )\times \mathscr A$ be  the configuration space and let $\cal M$ be the
subspace of critical points of the Seiberg-Witten functional. We
define $\Lambda :=\{ SW(\varphi_{\infty}, A_{\infty}):
(\varphi_{\infty}, A_{\infty})\in {\cal M}\}$.   By the
compactness result in \cite{Jost.Variational} and Lemma \ref{Lojasiewicz}, we
know that $\Lambda$ is discrete. For each $\lambda \in \Lambda$,
let $\cal M_{\lambda}$ be the subset of critical points $(\varphi_\infty,A_\infty)$ with $\mathcal{SW}(\varphi_\infty,A_\infty)= \lambda$,
 and ${\cal A}_{\lambda}$ the subset of $\cal A$
such that $SW(\varphi (t), A(t))\to \lambda$. Then ${\cal
A}=\cup_{\lambda \in{\Lambda}} {\cal A}_{\lambda}$ and ${\cal
M}=\cup_{\lambda \in \Lambda } {\cal M}_{\lambda}$.  As a
consequence of Theorem 2, the Seiberg-Witten flow defines a
continuous $\mathscr G$-equivariant flow. Furthermore, the Seiberg-Witten
flow defines a deformation retraction $\Phi :[0,\infty ]\times
{\cal A}_{\lambda}\to {\cal A}_{\lambda}$ of ${\cal A}_{\lambda}$
onto $\cal M_{\lambda}$.

\medskip

It is a very interesting question when the  unique limit
$(\varphi_{\infty}, A_{\infty})$  of the Seiberg-Witten flow for
some initial data is a Seiberg-Witten monopole. By Lemma \ref{uniqueness}, if
the initial data $(\varphi_0, A_0)$ is sufficiently close to a
non-trivial Seiberg-Witten monopole, the flow will converge to a non-trivial Seiberg-Witten monopole which is close to the original non-trivial monopole. If the scalar curvature $S$ is everywhere
non-negative, the Seiberg-Witten equations (\ref{sweq}) admit only
the trivial solutions $\varphi =0$ and $F_A^+=0$, and equations
(\ref{eq1})-(\ref{eq2}) admit only trivial-type solutions with $\varphi =0$.
Thus, the flow can only converge to a trivial critical point.

\medskip

The paper is organized as follows: In Section 2, we establish some
preliminary estimates. In Section  3, we show the local existence
of the flow. In Section 4, we show global existence and complete
the proof of Theorem 1. In Section 5, we consider the limiting
behaviour of the flow and prove Theorem 2. Finally, in Section 6,
we present a brief note about analogous results for the flow of
the perturbed Seiberg-Witten functional.

\section{Preliminary estimates}

The familiar Sobolev spaces of functions on Euclidean spaces can
be extended to the geometrical context.  Given a connection $\nabla ^E :\Omega ^0 (E) \to \Omega ^1 (E)$ on a vector bundle $E$, we can extend it to the well-known exterior covariant derivative $d_A:\Omega ^p (E) \to \Omega ^{p + 1} (E)$. There is another extension of $\nabla ^E$, called the iterated covariant derivative
\[
\nabla : \otimes ^p T^* M \otimes E \to  \otimes ^{p + 1} T^* M \otimes E.
\]
We then define
\[
\left\| \varphi  \right\|_{W^{k,p}(M) }  = \left( {\sum\limits_{n
= 0}^k {\int_M {\left| {\nabla _{ref}^{(n)}  \varphi } \right|^p
\,dV} } } \right)^{\frac{1} {p}},
\]
where $\nabla _{ref}$ is a given reference connection and $\nabla
_{ref}^{(n)}$ denotes $n$ iterations of $\nabla _{ref}$ (we use the exponent $n$ without the brackets to denote the
$n^{th}$ component). It is a straightforward calculation to show
that different choices of reference connection lead to equivalent
norms. We define $\left\| A  \right\|_{W^{k,p}}$ similarly, where
the reference connection is simply the standard connection on
forms induced by the Levi-Civita connection. We define, as usual,
$H^k=W^{k,2}$. We also have the parabolic  spaces $L^p ([0,T];
W^{k,p}(M))$, which require that the function $t \to \left\|
{\varphi (t)} \right\|_{W^{k,p}}$ is in $L^p$ over $[0,T]$. In particular,
\[
\left\| \varphi  \right\|_{L^2([0,T]; L^2(M) )}^2  = \int_{M
\times \left[ {0,T} \right]} {\left| \varphi  \right|^2 }\,dVdt.
\]
We make use of another Weitzenb\"ock formula on $p$-forms (one that is distinct from (\ref{Weitzenbock})). We have the covariant
Laplacian $\nabla _M^* \nabla _M$ and the Hodge Laplacian
$\Delta=(dd^*+d^*d)$ (which has opposite sign to the standard
Laplace operator on $M$) . They are related by
\begin{equation} \label{weitzenbockforms}
\nabla _M^* \nabla _M \beta - \Delta \beta = R_M \#  \beta,
\end{equation}
where $\beta$ is any $p$-form, $R_M$ is the curvature of the
Levi-Civita connection, and \# represents some multilinear map
with smooth coefficients (so that importantly $\left| {R_M \#
\beta} \right| \leqslant c\left| {R_M } \right|\left|
\beta\right|$). See (\cite{Jost.Riemannian}) for details.

We first establish a bound on $\left| \varphi  \right|$. Let
$S_{0}=\min \left\{S(x):x\in M\right\}$. Of course, if $S_0 \geq 0$, the Seiberg-Witten equations admit only the trivial
(reducible) solutions $\varphi=0$ and $F^{+}_{A}=0$.
\begin{Lemma1} \label{phibound}
Let $(\varphi,A)$ be a solution of (\ref{flow1})-(\ref{flow2})
on $M \times [0,T)$, and write $m=\mathop {\sup }\limits_{x \in M}
\left| {\varphi _0 } \right|$. Then for all $t\in [0,T)$, we have
\begin{equation} \label{varphibound}
\mathop {\sup_{x\in M} } \left| {\varphi} (x,t) \right| \leqslant
\max \{ m,\sqrt {\left| {S_0 } \right|} \}.
\end{equation}
\end{Lemma1}
\emph{Proof}. We note the following identity:
\begin{equation} \label{laplaceidentity}
\Delta \left| \varphi  \right|^2  = 2\operatorname{Re}
\left\langle {\nabla_A^* \nabla_A \varphi ,\varphi } \right\rangle
- 2\left| {\nabla_A \varphi } \right|^2,
\end{equation}
which holds for any metric connection $\nabla_A$ (see 3.2.7 of \cite{Jost.Riemannian}). Using this identity, we have
\begin{align*}
\frac{\partial } {{\partial t}}\left| \varphi  \right|^2 & =
2\operatorname{Re} \left\langle {\frac{{\partial \varphi }}
{{\partial t}},\varphi } \right\rangle\\
& = 2\operatorname{Re} \left\langle { - \nabla _A^* \nabla _A
\varphi  - \frac{1}
{4}\left[ {S + \left| \varphi  \right|^2 } \right]\varphi ,\varphi } \right\rangle\\
& =  - \Delta \left| \varphi  \right|^2  - 2\left| {\nabla _A
\varphi } \right|^2  - \frac{1} {2}\left[ {S + \left| \varphi
\right|^2 } \right]\left| \varphi  \right|^2.
\end{align*}
Let $b$ be any constant with $0<b<T$. Suppose $\phi (x,t)$ attains
its maximum point at $(x_0, t_0)\in M\times [0,b]$ such that $t_0$
is the first time the maximum is reached, i.e.
\[
|\varphi (x_0,t_0)|=\max_{x\in M, 0\leq t\leq b} |\varphi (x,t)|
\]

If $|\varphi (x_0,t_0)|\leq   \max \{m, \sqrt {\left| {S_0 } \right|}  \}$, the
claim is proved. Otherwise,
\[
|\varphi (x_0,t_0)|> \max \{m, \sqrt {\left| {S_0 } \right|}  \}.
\]
By the continuity of $\varphi$ on $M\times [0,b]$, there is a
parabolic cylinder $U\times [t_1, t_2]$ inside $M\times [0,b]$
with $t_1<t_0\leq t_2$ such that
\[
|\varphi (x,t)|\geq \max \{m, \sqrt {\left| {S_0 } \right|}  \}, \quad \forall
(x,t)\in U\times [t_1,t_2].
\]
Then for all $(x,t)\in U\times [t_1, t_2]$ we have
\[\frac{\partial } {{\partial t}}\left| \varphi  \right|^2  + \Delta \left| \varphi  \right|^2   \leq
0.
\]

By the strong parabolic maximum principle, $|\phi (x,t)|$ must be
a constant. This is impossible. $\square$ \vspace{5mm}

We have the following energy inequality.

\begin{Lemma2} \label{energyinequalitylemma}
Let $(\varphi,A)$ be a solution of (\ref{flow1})-(\ref{flow2})
on $M \times [0,T)$. Then
\begin{equation} \label{energyinequality}
\frac{d} {{dt}}\mathcal{SW}(\varphi (t),A(t)) =  - \int_M {\left[ 2{\left|
{\frac{{\partial \varphi }} {{\partial t}}} \right|^2  + \left|
{\frac{{\partial A}} {{\partial t}}} \right|^2 } \right]}
\leqslant 0.
\end{equation}

\end{Lemma2}
\emph{Proof}. For any $\psi$, we compute

\begin{align*}
&\quad \left. {\frac{d} {{d\varepsilon }}} \right|_{\varepsilon  =
0}
\mathcal{SW}(\varphi + \varepsilon \psi ,A) \\
&= 2\int_M {\left( {\operatorname{Re} \left\langle {\nabla _A^*
\nabla _A \varphi ,\psi } \right\rangle  + \frac{1} {4}\left[ {S +
\left| \varphi \right|^2 } \right]\operatorname{Re} \left\langle
{\varphi ,\psi } \right\rangle } \right)}\\
 &= 2\int_M
{\operatorname{Re} \left\langle {\nabla _A^* \nabla _A \varphi  +
\frac{1}{4}\left[ {S + \left| \varphi  \right|^2 } \right]\varphi
,\psi } \right\rangle },
\end{align*}
and for $B \in i\Lambda ^1 M$,
 \begin{align*}
 &\quad \left. {\frac{d} {{d\varepsilon }}} \right|_{\varepsilon  = 0}
\mathcal{SW}(\varphi ,A + \varepsilon B)\\
&= \left. {\frac{d} {{d\varepsilon }}} \right|_{\varepsilon  = 0}
\int_M {\left\langle {\nabla _{A + \varepsilon B} \; \varphi ,\nabla
_{A + \varepsilon B} \; \varphi } \right\rangle  + \left\langle {F_{A
+ \varepsilon B}^ + ,F_{A + \varepsilon B}^ +  } \right\rangle  +
\frac{S} {4}\left| \varphi \right|^2 + \frac{1} {8}\left| \varphi
\right|^4 }\\
&= \left. {\frac{d} {{d\varepsilon}}} \right|_{\varepsilon = 0}
\int_M {\left( {\left\langle {\nabla _A \varphi  +
\frac{1}{2}\varepsilon B\varphi ,\nabla _A \varphi  +
\frac{1}{2}\varepsilon B\varphi } \right\rangle  + \left\langle
{F_A^ +   + \varepsilon (dB)^ + ,F_A^ +   + \varepsilon (dB)^ + }
\right\rangle } \right)} \\
&  = 2\int_M {\left( {\frac{1}{2}\operatorname{Re} \left\langle
{\nabla _A \varphi ,B\varphi } \right\rangle  + \left\langle {F_A^
+  ,(dB)^ +  } \right\rangle } \right)}\\ & = 2\int_M {\left(
{\left\langle {\frac{i}{2}\operatorname{Im} \left\langle  {\nabla
_A \varphi ,\varphi } \right\rangle ,B} \right\rangle  +
\left\langle {F_A^+ ,dB} \right\rangle }
\right)},\\
& = 2\int_M {\left( {\left\langle {\frac{i}{2}\operatorname{Im}
\left\langle  {\nabla _A \varphi ,\varphi } \right\rangle
+\frac{1}{2}d^*F_A ,B} \right\rangle } \right)}, \end {align*}
where we have used that $d^* (dA)^ +   = \frac{1}{2}d^* dA$.
Noting that $\frac{{\partial \varphi }} {{\partial t}} = \psi$ and
$\frac{{\partial A}}{{\partial t}} = B$, the result follows.
$\square$ \vspace{5mm}

Next, integrating  (\ref{energyinequality})  in time gives
\begin{equation} \label{timeintegralbound}
\int_0^T  {\left[ 2{\left\| {\frac{{\partial \varphi }} {{\partial
t}}} \right\|_{L^2 }^2  + \left\| {\frac{{\partial A}} {{\partial
t}}} \right\|_{L^2 }^2 } \right]}  =  \mathcal{SW}(\varphi_0,A_0)
- \mathcal{SW}( \varphi (T),A(T) ).
\end{equation}
That is,
\[
\frac{{\partial \varphi }} {{\partial t}},\frac{{\partial A}}
{{\partial t}} \in L^2 ([0,T]; L^2(M) ).
\]
From the Seiberg-Witten functional (\ref{swfunc}) we see that
\[
\left\| {\nabla _A \varphi } \right\|_{L^2 }^2  + \left\| {F_A^ +
} \right\|_{L^2 }^2  + \frac{1}{4}\int_M {S\left| \varphi
\right|^2 }  \leqslant \mathcal{SW}(\varphi ,A) \leqslant
\mathcal{SW}(\varphi _0 ,A_0 )
\]
\[
 \Rightarrow \left\| {\nabla _A \varphi } \right\|_{L^2 }^2  + \left\| {F_A^ +  } \right\|_{L^2 }^2  \leqslant \mathcal{SW}(\varphi _0 ,A_0 ) - \frac{1}
{4}\int_M {S\left| \varphi  \right|^2 }
\]
\[
 \Rightarrow \left\| {\nabla _A \varphi } \right\|_{L^2 }^2  + \left\| {F_A^+ } \right\|_{L^2 }^2  \leqslant c,
\]
since $S$ and ${\left| \varphi  \right|}$ are bounded. This
implies that
\[
\nabla _A \varphi  \in L^\infty  ([0,T]; L^2(M) ).
\]
Furthermore, since from (\ref{selfdualrelation}), $\left\| {F_A^ +  } \right\|_{L^2
}^2=\frac{1}{2}\left\| {F_A } \right\|_{L^2 }^2 +c$, we also have
\[
F_A  \in L^\infty ([0,T]; L^2(M) ).
\]

\section{Local Existence} \label{localexistsection}

In this section, we show the existence of a classical (smooth)
solution of the  system (\ref{flow1})-(\ref{flow2}) on $M \times [0,T)$ for some
$T>0$. What we would like to do is to make the system parabolic by adding the term $dd^* A$ to (\ref{flow2}), since
this would give us the Laplacian $\Delta A$. Note that $\Delta
=dd^*+d^*d$ denotes the Hodge Laplacian. Fortunately, this extra
term points along the gauge orbit of $A$ since it is the
derivative of a function on $M$.

In local coordinates, we write
\[d_{\tilde A}=d+\tilde A=d_{A_0}+\tilde a,\quad d_A=d_{A_0}+a.\]
Then we consider the following system of equations:
\begin{equation} \label{flow3}
\frac{{\partial \tilde\varphi }} {{\partial t}}  = -\nabla
_{\tilde A}^* \nabla _{\tilde A} \tilde\varphi- \frac{1} {4}\left[
{S + \left| {\tilde \varphi } \right|^2 } \right]\tilde \varphi  +
\frac{1}{2}d^* \tilde a \tilde \varphi,
\end{equation}
\begin{equation} \label{flow4}
\frac{{\partial \tilde a}} {{\partial t}}  = - d^*F_{\tilde A}
-i\operatorname{Im} \left\langle {\nabla _{\tilde A} {\tilde
\varphi} ,{\tilde \varphi} } \right\rangle -dd^*\tilde a,
\end{equation}
with initial value $\tilde a (0)= 0$ and $\tilde \varphi
(0)=\varphi_0$.

Since $F_{\tilde A}=F_{A_0}+d\tilde a$ and
 \[-\nabla
_{\tilde A}^* \nabla _{\tilde A} \tilde\varphi=-\nabla _{A_0}^*
\nabla _{A_0} \tilde\varphi+ \tilde a\#\nabla_{A_0}\tilde\varphi
+\tilde a\#\tilde a\tilde \varphi +\nabla_{A_0}\tilde
a\#\tilde\varphi ,\] the system (\ref{flow3}) and (\ref{flow4}) is
a quasilinear parabolic system. Thus by standard PDE theory
there is a local smooth solution $(\tilde\varphi ,\tilde a)$ on $M
\times [0,T)$ for some $T>0$, given smooth initial data. See for
instance \S III.4 of \cite{Eidel'man}. However, since $\tilde a$ is not
bounded, we do not yet have global existence for the system
(\ref{flow3}) and (\ref{flow4}).

We next claim that the system (\ref{flow3})-(\ref{flow4}) is
gauge equivalent to our original system (\ref{flow1})-(\ref{flow2}). By standard ODE theory, there is a local smooth
solution $f$ to  the equation
\[
g^{-1} \frac {d g}{dt}  = - d^* \tilde a,
\]
\[
g(0) = I,
\]
where $( \tilde \varphi,\tilde a)$ is a solution to (\ref{flow3})
and (\ref{flow4}). Since $d^* \tilde a\in \Omega^0({\mathcal
U}(1))$, it is easy to check that $g^{-1}$ and ${\tilde g}^t$
satisfies
\[ \frac {d G}{dt} =d^*\tilde a\, G, \quad G(0)=I.
\]
Therefore $g^{-1}={ \tilde g}^t$. Hence, $g$ is a gauge
transformation.

Note that locally on the manifold we can write $g=e^{if}$ for some
real-valued function $f$. Then $d(g^{-1}dg)=0$, and   $g$ also
satisfies the equation
\begin{equation} \label{gaugechange}
2\frac{\partial }{{\partial t}}\left( g^{-1}dg\right)=-dd^* \tilde a.
\end{equation}
 Given our local solution $(\tilde\varphi , \tilde a)$ to
(\ref{flow3})-(\ref{flow4}) on $M \times \left[ {0,T}
\right)$, we solve equation (\ref{gaugechange}) to obtain our
gauge transformation $g(t)$.
 Set \[(\tilde \varphi , d_{\tilde
A})=(g^*(\varphi ), g^*(d_A))=(g^{-1}\varphi , d_A-2g^{-1}dg ).\]
Applying this gauge transformation, we obtain a local solution
\[(\varphi ,d_A) = (g  \tilde\varphi ,d_{\tilde A }- 2g^{-1}dg)\] to
our original system (\ref{flow1}) and (\ref{flow2}) on $M \times
\left[ {0,T} \right)$ as shown below.

Note that $F_{\tilde
A}=F_A$ since $d(g^{-1}dg )=0$  and that \[ g  \nabla_{\tilde A}^*
\circ g^{-1} =\nabla^* _{A},\quad g\nabla _{\tilde A }\circ
g^{-1}=\nabla _{A},\quad |\varphi |=|\tilde \varphi |.
\]
Then we have
\[\operatorname{Im} \left\langle {\nabla
_{\tilde A}  \tilde \varphi , \tilde \varphi } \right\rangle
=\operatorname{Im} \left\langle {g \nabla _{\tilde A} \tilde
\varphi ,g  \tilde \varphi } \right\rangle =\operatorname{Im}
\left\langle {\nabla _A \varphi ,\varphi } \right\rangle .\]
 We compute
\begin{align*}
&\frac{\partial A}{\partial t} =\frac {\partial \tilde a}{\partial t}-2\frac {\partial }{\partial t}(g^{-1}dg)\\
& = - d^* F_{\tilde A}   - i\operatorname{Im} \left\langle {
\nabla _{ \tilde A} \tilde \varphi , \tilde\varphi } \right\rangle
-d d^* \tilde a-
2\frac{\partial } {{\partial t}}\left( g^{-1}dg \right) \\
& = - d^* F_{\tilde A}   - i\operatorname{Im} \left\langle {
\nabla
_{\tilde A} \tilde\varphi , \tilde \varphi } \right\rangle  \\
&= - d^* F_{A}  - i\operatorname{Im} \left\langle {\nabla _A
\varphi ,\varphi } \right\rangle
\end{align*}
and
\begin{align*}
&\frac{\partial \varphi}{\partial t}
=\frac{\partial}{{\partial t}}(g \tilde\varphi )=\frac{{\partial g }}{{\partial t}}  \tilde\varphi  +g \frac {\partial
\tilde \varphi}{\partial t} \\
& = \frac{{\partial g }}{{\partial t}} \tilde\varphi  + g  \left(
-\nabla _{\tilde A}^* \nabla_{\tilde A}\tilde\varphi- \frac{1}
{4}\left[ {S + \left| {\tilde \varphi } \right|^2 } \right] \tilde
\varphi  + \frac{1}{2}d^* \tilde a \tilde \varphi, \right)\\
& = - g  \nabla_{\tilde A}^*\circ  g^{-1}\circ  g\nabla _{\tilde A
}\circ g^{-1}g \tilde \varphi - \frac{1}
{4}\left[ S + \left| { \tilde\varphi } \right|^2  \right] g \tilde \varphi \\
& = -\nabla _{A}^* \nabla _{A}  \varphi  - \frac{1} {4}\left[ S +
\left| { \varphi } \right|^2  \right] \varphi .
\end{align*}
Thus, we have shown the existence of the local solution of (\ref{flow1})-(\ref{flow2}).

\begin{Local} \label{Local}
For any given smooth initial data $(\varphi_0, A_0)$, equations
(\ref{flow1}) and (\ref{flow2}) admit a unique local smooth
solution on $M \times \left[ 0,T \right)$ for some $T>0$.
\end{Local}

We suppose that $T$ is maximal, that is, the solution cannot be smoothly extended
beyond time $T$, and contradict this assumption in the next
section.

\section{Global Existence}

 In this section, we show that our
local solution can be extended to a global solution, without
restrictions on the manifold, bundles, or initial data. The
obstruction to extending the local solution of (\ref{flow1})-(\ref{flow2}) on $M \times [0,T)$ to a global solution on $M
\times [0,\infty )$ is that it may cease to be smooth in finite
time.   Throughout this section, $(\varphi,A)$ will represent our
smooth local solution to the flow on $M \times [0,T)$. For
notational simplicity, we adopt the convention that $c$ and its
variants denote positive constants, which can change from line to
line.

We next compute an estimate for $\frac{\partial } {{\partial
t}}\left( {\left| {\nabla _A \varphi } \right|^2  + \left| {F_A }
\right|^2 } \right)$.

\begin{Lemma3} \label{firstderivativeestimate}
There exist positive constants $c,c'$ such that the following
estimate holds:
\begin{align*}
&\frac{\partial } {{\partial t}}\left( {\left| {\nabla _A \varphi
} \right|^2  + \left| {F_A } \right|^2 } \right) + \Delta \left(
{\left| {\nabla _A \varphi } \right|^2  +\left| {F_A } \right|^2 }
\right)  \\
&\leqslant - c'\left( {\left| {\nabla _A^2 \varphi } \right|^2  +
\left| {\nabla F_A } \right|^2 } \right) + c\left( {\left| {F_A }
\right| + 1} \right)\left( {\left| {\nabla _A \varphi } \right|^2
+ \left| {F_A } \right|^2 +1} \right).
\end{align*}
\end{Lemma3}
\emph{Proof}. We first consider $\left| {\nabla _A \varphi }
\right|^2$.
\begin{align}
&\frac{\partial } {{\partial t}}\left| {\nabla _A \varphi }
\right|^2   = 2\operatorname{Re} \left\langle {\nabla _A
\frac{{\partial \varphi }} {{\partial t}} + \left( {\frac{\partial
} {{\partial t}}\nabla _A } \right)\varphi ,\nabla _A \varphi }
\right\rangle \nonumber\\
& \label{varphiderivatives}
 =  - 2\operatorname{Re} \left\langle {\nabla _A \nabla _A^* \nabla _A \varphi ,\nabla _A \varphi } \right\rangle
  - \frac{1}{2}\operatorname{Re} \left\langle {\nabla _A \left[ {S + \left| \varphi  \right|^2 } \right]\varphi ,\nabla _A \varphi }
  \right\rangle \nonumber\\
  &\quad + \operatorname{Re} \left\langle {\frac{{\partial A}}
{{\partial t}}\varphi ,\nabla _A \varphi } \right\rangle.
\end{align}
Recall that we denote the curvature of the induced connection on
$\mathcal{S^+}$ by $\Omega_A$ with $A=A_0+a$, $a\in i\Lambda^1 M$.
We have the well-known Ricci formula
\begin{equation} \label{Ricciformula}
\nabla _A^{(n)} \nabla _A^* \nabla _A \varphi  = \nabla _A^*
\nabla _A \nabla _A^{(n)} \varphi  + \sum\limits_{j + k = n}
{\left( {\nabla_M^{(j)} R_M  + \nabla _M^{(j)} \Omega _A }
\right)\# } \nabla _A^{(k)} \varphi,
\end{equation}
where $R_M$ represents the Riemannian curvature of $M$ (see e.g. 2.2 of \cite{Hong.Zheng}). Then
 \begin{align}
- 2\operatorname{Re} \left\langle {\nabla _A \nabla _A^* \nabla _A
\varphi ,\nabla _A \varphi } \right\rangle  \leqslant & -
2\operatorname{Re} \left\langle {\nabla _A^* \nabla _A \nabla _A
\varphi ,\nabla _A \varphi } \right\rangle +c\left| {F_A }
\right|\left| {\nabla _A \varphi } \right|^2\nonumber \\
&+ c\left| {\nabla
_M F_A } \right|\left| {\nabla _A \varphi } \right|
 \label{curvature} + c\left| {\nabla _A \varphi } \right|^2 +
c\left| {\nabla _A \varphi } \right|,
\end{align}
where we note that the nonconstant portion of $\Omega_A$ is
$F_A$. We deal with the first term in  (\ref{varphiderivatives})
by applying (\ref{laplaceidentity}) to $\nabla _A \varphi$:
\[
 - 2\operatorname{Re} \left\langle {\nabla _A^* \nabla _A \nabla _A \varphi ,\nabla _A \varphi } \right\rangle
  =  - \Delta \left| {\nabla _A \varphi } \right|^2  - 2\left| {\nabla _A^{(2)} \varphi } \right|^2.
\]
Considering now the second term in  (\ref{varphiderivatives}), we
note that by the metric compatibility we have
\[
d\left| \varphi  \right|^2  = \left\langle {\nabla _A \varphi
,\varphi } \right\rangle  + \left\langle {\varphi ,\nabla _A
\varphi } \right\rangle  = 2\operatorname{Re} \left\langle {\nabla
_A \varphi ,\varphi } \right\rangle
\]
and so
\begin{align*}
&  -\frac{1}{2}\operatorname{Re} \left\langle {\nabla _A \left[ {S + \left| \varphi  \right|^2 } \right]\varphi ,\nabla _A \varphi } \right\rangle \\
& = -\frac{1}{2}\operatorname{Re} \left\langle { \left[ {S + \left| \varphi  \right|^2 } \right]\nabla _A \varphi  +dS\varphi
+ d\left| \varphi  \right|^2 \varphi ,\nabla _A \varphi } \right\rangle \\
& = - \frac{1} {2}\left[ {S + \left| \varphi  \right|^2 }
\right]\left| {\nabla _A \varphi } \right|^2  - \frac{1}
{2}\operatorname{Re} \left\langle {dS\varphi ,\nabla _A \varphi } \right\rangle
- \operatorname{Re} \left\langle {\operatorname{Re} \left\langle {\nabla _A \varphi ,\varphi } \right\rangle \varphi,\nabla _A \varphi } \right\rangle \\
& =  - \frac{1} {2}\left[ {S + \left| \varphi  \right|^2 }
\right]\left| {\nabla _A \varphi } \right|^2  - \frac{1}
{2}\operatorname{Re} \left\langle {dS\varphi ,\nabla _A \varphi } \right\rangle
- \left| {\operatorname{Re} \left\langle {\nabla _A \varphi ,\varphi } \right\rangle } \right|^2 \\
& \leqslant c\left| {\nabla _A \varphi } \right|^2  + c\left|
{\nabla _A \varphi } \right|,
\end{align*}
where we have used that
\begin{align*}
& - \operatorname{Re} \left\langle {\operatorname{Re} \left\langle {\nabla _A \varphi ,\varphi } \right\rangle \varphi ,\nabla _A \varphi } \right\rangle
= - \operatorname{Re} \left\langle {\operatorname{Re} \left\langle {\nabla _A^j \varphi ,\varphi } \right\rangle \varphi ,\nabla _A^j \varphi } \right\rangle\\
& = - \operatorname{Re} \left\langle {\nabla _A^j \varphi ,\varphi } \right\rangle \operatorname{Re} \left\langle {\varphi ,\nabla _A^j \varphi } \right\rangle
 = -\left| {\operatorname{Re} \left\langle {\varphi ,\nabla _A \varphi } \right\rangle } \right|^2.\\
\end{align*}
Finally, for the third term in (\ref{varphiderivatives}),
\begin{align*}
\operatorname{Re} \left\langle {\frac{{\partial A}}{{\partial
t}}\varphi ,\nabla _A \varphi } \right\rangle & = -
\operatorname{Re} \left\langle {d^* F_A \varphi ,\nabla _A \varphi
} \right\rangle
 - \operatorname{Re} \left\langle {i\operatorname{Im} \left\langle {\nabla _A \varphi ,\varphi } \right\rangle \varphi ,\nabla _A \varphi } \right\rangle\\
& \leqslant c\left| {\nabla _M F_A } \right|\left| {\nabla _A \varphi } \right| + c\left| {\nabla _A \varphi } \right|^2. \\
\end{align*}
Combining all of the above we ultimately find
\begin {align}
 \frac{\partial } {{\partial t}}\left| {\nabla _A
\varphi } \right|^2  \leqslant &  - \Delta \left| {\nabla _A
\varphi } \right|^2  - \left| {\nabla _A^2 \varphi } \right|^2 +
c\left| {\nabla _M F_A } \right|\left| {\nabla _A \varphi }
\right|\nonumber \\
&+ c\left| {F_A } \right|\left| {\nabla _A \varphi } \right|^2 +
c\left| {\nabla _A \varphi } \right|^2 + c\left| {\nabla _A
\varphi } \right|.\label{estimate1}
\end{align}
We next consider $\left| {F_{A} } \right|^2$.
\begin{align*}
&\frac{\partial } {{\partial t}}\left| {F_A } \right|^2  =
\frac{\partial } {{\partial t}}\left| {dA} \right|^2  =
2\left\langle {d\frac{{\partial A}} {{\partial t}},dA}
\right\rangle
\\
&
 = 2\left\langle {d\left[ { - d^* F_A  - i\operatorname{Im} \left\langle {\nabla _A \varphi ,\varphi } \right\rangle } \right],F_A } \right\rangle
\\
&
 = 2\left\langle { - \Delta F_A  - id\operatorname{Im} \left\langle {\nabla _A \varphi ,\varphi } \right\rangle ,F_A } \right\rangle
\\
&
 =  - 2\left\langle {\Delta F_A ,F_A } \right\rangle  - 2\left\langle {id\operatorname{Im} \left\langle {\nabla _A \varphi ,\varphi } \right\rangle ,F_A } \right\rangle,
\end{align*}
where we have utilized the Bianchi identity $dF_A  = 0$, giving
$dd^* F_A  = \Delta F_A$. Applying the Weitzenb\"ock formula
(\ref{weitzenbockforms}) and recalling (\ref{laplaceidentity}),
\begin{align*}
2\left\langle {\Delta F_A ,F_A } \right\rangle  &\leqslant
2\left\langle {\nabla _M^* \nabla _M F_A ,F_A } \right\rangle +
c\left| {F_A } \right|^2
\\
& =  - \Delta \left| {F_A } \right|^2  - 2\left| {\nabla _M F_A } \right|^2 + c\left| {F_A } \right|^2.
\end{align*}
Then using metric compatibility
\begin{align} \label{Fmetric}
&\quad d\operatorname{Im} \left\langle {\nabla _A \varphi ,\varphi
} \right\rangle
 = d\left( {\operatorname{Im} \left\langle {\nabla _A^j \varphi ,\varphi } \right\rangle dx^j } \right)\nonumber \\
& = d_k (\operatorname{Im} \left\langle {\nabla _A^j \varphi ,\varphi } \right\rangle )dx^k  \wedge dx^j\nonumber\\
& = \sum\limits_{k > j} {\left( {d_k (\operatorname{Im} \left\langle {\nabla _A^j \varphi ,\varphi } \right\rangle )
- d_j (\operatorname{Im} \left\langle {\nabla _A^k \varphi ,\varphi } \right\rangle )} \right)dx^k  \wedge dx^j }\nonumber\\
& = \sum\limits_{k > j} { ( {\operatorname{Im} \left\langle {(\nabla _A^k \nabla _A^j  - \nabla _A^j \nabla _A^k )\varphi
,\varphi } \right\rangle + \operatorname{Im} \left\langle {\nabla_A^j \varphi ,\nabla _A^k \varphi } \right\rangle - \operatorname{Im} \left\langle {\nabla _A^k \varphi ,\nabla _A^j \varphi } \right\rangle }  )dx^k  \wedge dx^j } \nonumber \\
& = \sum\limits_{k > j} {\left( {\operatorname{Im} \left\langle {\Omega_A^{kj} \varphi ,\varphi } \right\rangle
+ 2\operatorname{Im} \left\langle {\nabla _A^j \varphi ,\nabla _A^k \varphi } \right\rangle } \right)dx^k  \wedge dx^j },
\end{align}
so that
\[
2\left\langle {id\operatorname{Im} \left\langle {\nabla _A \varphi
,\varphi } \right\rangle ,F_A } \right\rangle  \leqslant c\left|
{F_A } \right|\left| {\nabla _A \varphi } \right|^2  + c\left|
{F_A } \right|^2+ c\left| {F_A } \right|.
\]
Finally, we have
\begin{equation} \label{estimate2}
\frac{\partial } {{\partial t}}\left| {F_A } \right|^2  \leqslant
- \Delta \left| {F_A } \right|^2  - 2\left| {\nabla F_A }
\right|^2 + c\left| {F_A } \right|\left| {\nabla _A \varphi }
\right|^2 + c\left| {F_A } \right|^2+ c\left| {F_A } \right|.
\end{equation}
We now combine (\ref{estimate1}) and (\ref{estimate2}):
\begin{align}
&\frac{\partial } {{\partial t}}\left( {\left| {\nabla _A \varphi
}
\right|^2  + \left| {F_A } \right|^2 } \right)\nonumber \\
&\leqslant
 - \Delta \left( {\left| {\nabla _A \varphi } \right|^2  + \left| {F_A } \right|^2 } \right)
 - 2\left( {\left| {\nabla _A^{(2)} \varphi } \right|^2  + \left| {\nabla _M F_A } \right|^2 }
 \right)\nonumber
\\& \label{combined} + c\left| {\nabla _M F_A } \right|\left|
{\nabla _A \varphi } \right| + c\left( {\left| {F_A } \right| + 1}
\right)\left( {\left| {\nabla _A \varphi } \right|^2  + \left|
{F_A } \right|^2 } \right) + c
\end{align}
where the first powers of $\left| {F_A } \right|$ and $\left|
{\nabla _A \varphi } \right|$ can be incorporated into a constant
since if they are larger than one, they are bounded by the second powers. We next have to deal with the derivatives of the curvature
that appear in (\ref{combined}). Fortunately, they can be
controlled by the term $ - 2\left| {\nabla_M F_A } \right|^2$
using Young's inequality
\[
 \left| {\nabla _M F_A } \right|\left| {\nabla _A \varphi } \right| \leqslant \frac{1}{2} \varepsilon \left| {\nabla _M F_A } \right|^2  + \frac{1}
{2\varepsilon }\left| {\nabla _A \varphi } \right|^2.
\]
Then if we choose $\varepsilon$ sufficiently small we have the
desired result. $\square$ \vspace{5mm}

Using local coordinates, let \[P_R (y,s) = \{ (x,t) \in M \times
(0,T):\left| {x - y} \right| < R, \quad s-R^2<t< s \}\] be a
parabolic cylinder of radius $R$ centered at  $(y,s)$.

\begin{Lemma4} \label{regularitytheorem}
Suppose $(\varphi, A)\in C^{\infty}( P_R (y,s))$ satisfies
(\ref{flow1})-(\ref{flow2}). Then there exist  constants
$\delta$ and $R_0$ such that if $R \leqslant R_0$ and
\[
\mathop {\sup }\limits_{0 <t < s } \int_{B_R (y)} {\left(
{\left| {\nabla _A \varphi } \right|^2  + \left| {F_A } \right|^2
} \right)\,dV} < \delta,
\]
then
\[
\mathop {\sup }\limits_{P_{R/2} (y,s)} \left( {\left| {\nabla _A
\varphi } \right|^2  + \left| {F_A } \right|^2 } \right) \leqslant
256R^{ - 4}.
\]
\end{Lemma4}
\emph{Proof.} The proof is similar to one in \cite {Hong.Tian},
but there are some differences. For completeness, we give details
here. We begin by choosing $r_0<R$ so that
\begin{equation} \label{r0def}
(R - r_0 )^4 \mathop {\sup }\limits_{P_{r_0 } (y,s)} \left(
{\left| {\nabla _A \varphi } \right|^2  + \left| {F_A } \right|^2
} \right) = \mathop {\max }\limits_{0 \leqslant r \leqslant R}
\left[ {(R - r)^4 \mathop {\sup }\limits_{P_r (y,s)} \left(
{\left| {\nabla _A \varphi } \right|^2  + \left| {F_A } \right|^2
} \right)} \right].
\end{equation}
Let
\[
e_0  = \mathop {\sup }\limits_{P_{r_0 } (y,s)} \left( {\left|
{\nabla _A \varphi } \right|^2  + \left| {F_A } \right|^2 }
\right) = \left( {\left| {\nabla _A \varphi } \right|^2  + \left|
{F_A } \right|^2 } \right)(x_0 ,t_0 )
\]
for some $(x_0 ,t_0 ) \in \bar P_{r_0 } (y,s)$. We claim that
\begin{equation} \label{eclaim}
e_0  \leqslant 16(R - r_0 )^{ - 4}.
\end{equation}
Then
\begin{align*}
&(R - r)^4 \mathop {\sup }\limits_{P_r (y,s)} \left( {\left|
{\nabla _A \varphi } \right|^2  + \left| {F_A } \right|^2 }
\right) \leqslant (R - r_0 )^4 \mathop {\sup }\limits_{P_{r_0 }
(y,s)} \left( {\left| {\nabla _A \varphi } \right|^2  + \left|
{F_A } \right|^2 } \right)
\\
&
 \leqslant 16(R - r_0 )^4 (R - r_0 )^{ - 4}  = 16
\end{align*} for any $r<R$. Choosing $r = \frac{1}{2}R$ in the
above, we have the required result.

We now prove (\ref{eclaim}). Define $\rho _0  = e_0^{ - 1/4}$ and
suppose by contradiction that $\rho _0  \leqslant \frac{1}{2}(R -
r_0 )$. We rescale variables via $x=x_0+\rho_0\tilde x$ and
$t=t_0+\rho_0^2 \tilde t$ and set
\[
\psi (\tilde x,\tilde t) = \varphi (x_0  + \rho _0 \tilde x,t_0  +
\rho _0^2 \tilde t),
\]
\[
B(\tilde x,\tilde t) = \rho _0 A(x_0  + \rho _0 \tilde x,t_0  +
\rho _0^2 \tilde t),
\]
giving
\[
\left| {\nabla _B \psi } \right|^2  = \rho _0^2 \left| {\nabla _A
\varphi } \right|^2,
\]
\[
\left| {F_B } \right|^2  = \rho _0^4 \left| {F_A } \right|^2.
\]
We define
\[
e_{\rho _0 } (\tilde x,\tilde t) = \left| {F_B } \right|^2  + \rho
_0^2 \left| {\nabla _B \psi } \right|^2  = \rho _0^4 \left(
{\left| {\nabla _A \varphi } \right|^2  + \left| {F_A } \right|^2
} \right)
\]
so that
\[
e_{\rho _0 } (\tilde x,\tilde t) \leqslant e_{\rho _0 } (0,0) = 1.
\]
We compute
\begin{align*}
\mathop {\sup }\limits_{\tilde P_1 (0,0)} e_{\rho _0 } (\tilde
x,\tilde t)
& = \rho _0^4 \mathop {\sup }\limits_{P_{\rho _0 } (x_0 ,t_0 )} \left( {\left| {\nabla _A \varphi } \right|^2  + \left| {F_A } \right|^2 } \right) \\
& \leqslant \rho _0^4 \mathop {\sup }\limits_{P_{\frac{{R + r_0 }}
{2}} (y,s)} \left( {\left| {\nabla _A \varphi } \right|^2  + \left| {F_A } \right|^2 } \right)\\
& = \rho _0^4 \left( {\frac{{R - r_0 }} {2}} \right)^{ - 4} \left(
{R - \frac{{R + r_0 }} {2}} \right)^4 \mathop {\sup
}\limits_{P_{\frac{{R + r_0 }}
{2}} (y,s)} \left( {\left| {\nabla _A \varphi } \right|^2  + \left| {F_A } \right|^2 } \right)\\
& \leqslant \rho _0^4 \left( {\frac{{R - r_0 }}
{2}} \right)^{ - 4} \left( {R - r_0 } \right)^4 e_0  = 16, \\
\end{align*}
where we have used that $P_{\rho _0 } (x_0 ,t_0 ) \subset
P_{\frac{{R + r_0 }} {2}} (y,s)$, and to get to the last line we
have used (\ref{r0def}). This implies that
\[
e_{\rho _0 }  = \rho _0^4 \left( {\left| {\nabla _A \varphi }
\right|^2  + \left| {F_A } \right|^2 } \right) \leqslant 16
\]
on $\bar P_1(0,0)$. By Lemma \ref{firstderivativeestimate},
\begin{align*}
(\frac{\partial } {{\partial t}}+\Delta )  \left( {\left| {\nabla
_A \varphi } \right|^2  + \left| {F_A } \right|^2 }  +1 \right)
\leqslant  c\left( {\left| {F_A } \right| + 1} \right) \left( {
{\left| {\nabla _A \varphi } \right|^2  + \left| {F_A } \right|^2
}   + 1} \right).
\end{align*}
Then
\begin{align*}
\left( {\frac{\partial }{{\partial \tilde t}} + \tilde \Delta }
\right) (e_{_{\rho _0 } }+\rho_0^4)
& = \rho _0^6 \left( {\frac{\partial }{{\partial t}} + \Delta } \right)\left( {\left| {\nabla _A \varphi } \right|^2  + \left| {F_A } \right|^2 } \right)\\
& \leqslant c\rho _0^6 \left( {\left| {F_A } \right| + 1}
\right)\left( {\left| {\nabla _A \varphi } \right|^2  + \left|
{F_A } \right|^2 }+1 \right)
\end{align*}
on $\bar P_1(0,0)$. Note that by assumption $\rho _0  < R$, and thus $\rho_0^2 |F_A|$ is bounded by a constant. Then
\[
\left( {\frac{\partial } {{\partial \tilde t}} + \tilde \Delta }
\right)\left( {e_{\rho _0 }  +  \rho_0^4} \right) \leqslant
c\left( {e_{\rho _0 }  + \rho_0^4} \right)
\]
for a constant $c>0$.  We apply Moser's Harnack inequality to give
\begin{align*}
1 + \rho_0^4 = e_{\rho _0 } (0,0) + \rho_0^4
& \leqslant c\int_{\tilde P_1 (0,0)} {e_{\rho _0 } d\tilde xd\tilde t} +c\rho_0^4 \\
& = c\rho _0^{ - 2} \int_{P_{\rho _0 } (x_0 ,t_0 )} {\left( {\left| {\nabla _A \varphi } \right|^2  + \left| {F_A } \right|^2 } \right)} dxdt  +c\rho_0^4\\
& \leqslant c\mathop {\sup }\limits_{0  \leqslant t \leqslant s }
\int_{B_R (y)} {\left( {\left| {\nabla _A \varphi } \right|^2  +
\left| {F_A } \right|^2 } \right)}
 +  cR^4 \\
& < c\delta +cR^4,
\end{align*}
where  we have used that $\rho _0  < R$. Now if we choose $R_0$
 and $\delta $  sufficiently small, we have
the desired contradiction. $\square$ \vspace{5mm}

\begin{Lemma5} \label{localenergyestimate}
Let $(\varphi , A)$ be a solution to (\ref{flow1})-(\ref{flow2}).
Writing
\[
\mathcal{SW}_{B_R (x_0 )} (\varphi ,A) = \int_{B_R (x_0 )} {\left|
{\nabla _A \varphi } \right|^2  + \frac{1} {2}\left| {F_A }
\right|^2  + \frac{S} {4}\left| \varphi  \right|^2  + \frac{1}
{8}\left| \varphi  \right|^4 },
\]
we have for any $x_0 \in M$ and ball of radius $R$,
\[
\mathop {\sup }\limits_{t_1 \leqslant t \leqslant t_2}
\mathcal{SW}_{B_R (x_0 )} (\varphi ,A) \leqslant
\mathcal{SW}_{B_{2R} (x_0 )} (\varphi(t_1) ,A(t_1) ) +
C_1(t_2-t_1)R^{ - 2},
\]
where $C_1$ is a constant.
\end{Lemma5}
\emph{Proof.} Let $\phi$ be a smooth test function with $\phi
\equiv 1$ on $B_R(x_0)$ and zero outside of $B_{2R}(x_0)$. We can
choose $\phi$ so that $0 \leqslant \phi  \leqslant 1$ and $\left|
{d\phi } \right| \leqslant cR^{ - 1}$. We compute
\begin{align*}
&\frac{1}{2}\frac{d}{{dt}}\int_M {\phi ^2 \left| {F_A } \right|^2}
= \int_M {\left\langle {\phi ^2 d\frac{{\partial A}}
{{\partial t}},F_A } \right\rangle }  \\
& = \int_M {\left\langle {\phi ^2 \frac{{\partial A}} {{\partial
t}},d^* F_A } \right\rangle }  - \int_M {\left\langle {d\phi ^2
\wedge \frac{{\partial A}}
{{\partial t}},F_A } \right\rangle } \\
& \leqslant \int_M {\left\langle {\phi ^2 \frac{{\partial A}}
{{\partial t}},d^* F_A } \right\rangle }  + 2\int_M {\phi \left|
{d\phi } \right|} \left| {\frac{{\partial A}}
{{\partial t}}} \right|\left| {F_A } \right| \\
& \leqslant \int_M {\left\langle {\phi ^2 \frac{{\partial A}}
{{\partial t}},d^* F_A } \right\rangle }  + \int_M \phi  ^2 \left|
{\frac{{\partial A}}
{{\partial t}}} \right|^2  + \int_M {\left| {d\phi } \right|} ^2 \left| {F_A } \right|^2, \\
\end{align*}
and similarly
\begin{align*}
\frac{d} {{dt}}\int_M {\phi ^2 \left| {\nabla _A \varphi }
\right|^2}  &=  2\int_M {\phi ^2 \operatorname{Re} \left\langle
{\nabla _A \frac{{\partial \varphi }} {{\partial t}},\nabla _A
\varphi } \right\rangle }  + \int_M {\phi ^2 \operatorname{Re}
\left\langle {\frac{{\partial A}} {{\partial t}}\varphi ,\nabla _A
\varphi } \right\rangle }
\\&
 = 2\int_M {\phi ^2 \operatorname{Re} \left\langle {\frac{{\partial \varphi }}
{{\partial t}},\nabla _A^* \nabla _A \varphi } \right\rangle }  -
2\int_M {\operatorname{Re} \left\langle {d\phi ^2  \otimes
\frac{{\partial \varphi }} {{\partial t}},\nabla _A \varphi }
\right\rangle } \\
&\quad  + \int_M {\phi ^2 \operatorname{Re} \left\langle
{\frac{{\partial A}} {{\partial t}}\varphi ,\nabla _A \varphi }
\right\rangle }. \end{align*}

Furthermore, in the above
\begin{align*}
-2\int_M {\operatorname{Re} \left\langle {d\phi ^2  \otimes
\frac{{\partial \varphi }} {{\partial t}},\nabla _A \varphi }
\right\rangle } & \leqslant 4\int_M {\phi \left| {d\phi }
\right|\left| {\frac{{\partial \varphi }}
{{\partial t}}} \right|} \left| {\nabla _A \varphi } \right|\\
& \leqslant 2\int_M {\phi ^2 \left| {\frac{{\partial \varphi }}
{{\partial t}}} \right|} ^2  + 2\int_M {\left| {d\phi } \right|^2 } \left| {\nabla _A \varphi } \right|^2\\
\end{align*}
and
\begin{align*}
&2\int_M {\phi ^2 \operatorname{Re} \left\langle {\frac{{\partial
\varphi }} {{\partial t}},\nabla _A^* \nabla _A \varphi }
\right\rangle } \\
& =  - 2\int_M {\phi ^2 \left| {\frac{{\partial
\varphi }}{{\partial t}}} \right|^2 }  - 2\int_M {\phi ^2
\operatorname{Re} \left\langle {\frac{{\partial \varphi
}}{{\partial t}},\frac{1}
{4}\left[ {S + \left| \varphi  \right|^2 } \right]\varphi } \right\rangle } \\
& =  - 2\int_M {\phi ^2 \left| {\frac{{\partial \varphi }}
{{\partial t}}} \right|^2 }  - \frac{d} {{dt}}\int_M {\phi ^2
\left[ {\frac{S} {4}\left| \varphi  \right|^2  + \frac{1}
{8}\left| \varphi  \right|^4 } \right]}. \\
\end{align*}
Thus
\begin{align*}
\frac{d} {{dt}}\int_M {\phi ^2 \left| {\nabla _A \varphi }
\right|^2 \leqslant } &   - \frac{d} {{dt}}\int_M {\phi ^2 \left[
{\frac{S} {4}\left| \varphi  \right|^2  + \frac{1} {8}\left|
\varphi
\right|^4 } \right]}  \\
&+ 2\int_M {\left| {d\phi } \right|^2 } \left| {\nabla _A \varphi
} \right|^2  + \int_M {\phi ^2 \operatorname{Re} \left\langle
{\frac{{\partial A}} {{\partial t}}\varphi ,\nabla _A \varphi }
\right\rangle }. \end{align*}
We next note that
\begin{align} \label{imaginary}
\phi ^2 \operatorname{Re} \left\langle {\frac{{\partial A}}
{{\partial t}}\varphi ,\nabla _A \varphi } \right\rangle
& = \phi
^2 \operatorname{Im} \frac{{\partial A_k }} {{\partial
t}}\operatorname{Im} \left\langle {\nabla _A^k \varphi ,\varphi }
\right\rangle  \nonumber\\
& = \phi ^2 \left\langle {\frac{{\partial A}}
{{\partial t}},i\operatorname{Im} \left\langle {\nabla _A \varphi
,\varphi } \right\rangle } \right\rangle
\end{align}
so that
\[
\phi ^2 \operatorname{Re} \left\langle {\frac{{\partial A}}
{{\partial t}}\varphi ,\nabla _A \varphi } \right\rangle  + \int_M
{\left\langle {\phi ^2 \frac{{\partial A}} {{\partial t}},d^* F_A
} \right\rangle }  =  - \int_M \phi  ^2 \left| {\frac{{\partial
A}} {{\partial t}}} \right|^2.
\]
From all of the above, we finally have
\begin{align*}
&\frac{d} {{dt}}\int_M {\phi ^2 \left( {\left| {\nabla _A \varphi
} \right|^2  + \frac{1} {2}\left| {F_A } \right|^2  + \frac{S}
{4}\left| \varphi  \right|^2  + \frac{1} {8}\left| \varphi
\right|^4 } \right)} \\
& \leqslant cR^{ - 2} \int_M {\left( {\left| {\nabla _A \varphi }
\right|^2  + \frac{1} {2}\left| {F_A } \right|^2 } \right)}.
\end{align*}

The result follows by integrating on $[t_1,t]$ and taking the
supremum over $t_1 \leqslant t \leqslant t_2$. $\square$
\vspace{5mm}

\begin{Lemma6} \label{higherderivativeestimate} Let $(\varphi, A)$ be a solution to (\ref{flow1})-(\ref{flow2})  in $M\times [0,T)$
with initial values $(\varphi_0 ,A_0)$. Suppose $\left| {\nabla _A
\varphi } \right| \leqslant K_1$ and $\left| F_A   \right|
\leqslant K_1$ in $M\times [0, T)$ for   some constant $K_1>0$.
Then for any positive integer $n\geq 1$, there is a constant
$K_{n+1}$ independent of $T$  such that
\[\left|
{\nabla _A^{(n+1)} \varphi } \right| \leqslant K_{n+1},\quad
\left| {\nabla _M^{(n)} F_A } \right| \leqslant K_{n+1}\quad \mbox
{in } M\times [0, T).\]
\end{Lemma6}
\emph{Proof.} We prove Lemma \ref{higherderivativeestimate} by
induction. We first claim that
\begin{align}
&\frac{\partial } {{\partial t}}\left( {\left| {\nabla _A^{(k+1)}
\varphi } \right|^2  + \left| {\nabla _M^{(k)} F_A } \right|^2 }
\right)+ c'_k\left( {\left| {\nabla _A^{(k+2)} \varphi } \right|^2
+ \left| {\nabla_M^{(k+1)} F_A } \right|^2 } \right)\nonumber
\\
  \leqslant & -\Delta \left( {\left| {\nabla _A^{(k+1)} \varphi }
\right|^2  + \left| {\nabla _M^{(k)} F_A } \right|^2 } \right)
\label{higherorder} + c_k\left( {\left| {\nabla _A^{(k+1)} \varphi
} \right|^2  + \left| {\nabla _M^{(k)} F_A } \right|^2 +1}
\right)
\end{align}
for all non-negative integers $k=0, 1, 2, 3,\cdots$.

From Lemma \ref{firstderivativeestimate} with the assumption of Lemma
\ref{higherderivativeestimate}, (\ref {higherorder}) holds for
$k=0$. Now, assume that (\ref{higherorder}) is true for $k=n-1$
and  $\left| {\nabla _A^{(k+1)} \varphi } \right| \leqslant
K_{k+1}$ and $\left| {\nabla _M^{(k)} F_A } \right| \leqslant
K_{k+1}$ for non-negative integers $k\leq n-1$. Then we will show
(\ref{higherorder}) and Lemma \ref{higherderivativeestimate} are also true for
all $n$.

From (\ref{flow1}), we have
\begin{align}
\frac{\partial } {{\partial t}}\left| {\nabla _A^{(n+1)} \varphi }
\right|^2  &= 2\operatorname{Re} \left\langle {\frac{\partial }
{{\partial t}}\left( {\nabla _A^{(n+1)} \varphi } \right),\nabla
_A^{(n+1)} \varphi } \right\rangle \nonumber
\\
&
 =- 2\operatorname{Re} \left\langle {\nabla _A^{(n+1)} \nabla _A^* \nabla _A \varphi ,\nabla _A^{(n+1)} \varphi } \right\rangle \nonumber \\
 & - \frac{1}
{2}\operatorname{Re} \left\langle {\nabla _A^{(n+1)} \left[ {S +
\left| \varphi  \right|^2 } \right]\varphi ,\nabla _A^{(n+1)}
\varphi } \right\rangle\nonumber
\\
&\label{higherordereq1}
 + 2\operatorname{Re} \left\langle {\left( {\frac{\partial }
{{\partial t}}\nabla _A^{(n+1)} } \right)\varphi ,\nabla
_A^{(n+1)} \varphi } \right\rangle.
\end{align}
From the Ricci formula (\ref{Ricciformula}) we have
\begin{align*}
&- 2\operatorname{Re} \left\langle {\nabla _A^{(n+1)} \nabla _A^*
\nabla _A \varphi ,\nabla _A^{(n+1)} \varphi } \right\rangle
\\
&\leqslant
 - 2\operatorname{Re} \left\langle {\nabla _A^* \nabla _A \nabla _A^{(n+1)} \varphi ,\nabla _A^{(n+1)} \varphi } \right\rangle
\\&
+ c\left| {\nabla _M^{(n+1)} F_A } \right|\left| {\nabla
_A^{(n+1)} \varphi } \right| + c\left| {\nabla _A^{(n+1)} \varphi
} \right|^2 + c\left| {\nabla _A^{(n+1)} \varphi } \right|,
\end{align*}
where we recall that the non-constant part of $\Omega _A$ is
$F_A$. From (\ref{laplaceidentity}),
\[
 - 2\operatorname{Re} \left\langle {\nabla _A^* \nabla _A \nabla _A^{(n+1)} \varphi ,\nabla _A^{(n+1)} \varphi } \right\rangle  =
  - \Delta \left| {\nabla _A^{(n+1)} \varphi } \right|^2  - 2\left| {\nabla _A^{(n + 2)} \varphi } \right|^2.
\]
Next, applying metric compatibility $n+1$ times, we find that the
${(n+1)}^{th}$ order term of $\partial ^{(n+1)} \left| \varphi
\right|^2$ is $2\operatorname{Re} \left\langle {\nabla _{_A
}^{(n+1)} \varphi ,\varphi } \right\rangle$ and
\[
\operatorname{Re} \left\langle { - \nabla _A^{(n+1)} \frac{1}
{4}\left[ {R + \left| \varphi  \right|^2 } \right]\varphi ,\nabla
_A^{(n+1)} \varphi } \right\rangle  \leqslant c\left| {\nabla
_A^{(n+1)} \varphi } \right|^2  + c\left| {\nabla _A^{(n+1)}
\varphi } \right|.
\]
For the final term in (\ref{higherordereq1}), noting that
$\frac{\partial }{{\partial t}}\nabla _A  = \frac{1}
{2}\frac{{\partial A}}{{\partial t}}$ involves derivatives of
$F_A$ and $\nabla_A \varphi$ and utilizing the product rule we
find
\begin{align*}
&2\operatorname{Re} \left\langle {\left( {\frac{\partial }
{{\partial t}}\nabla _A^{(n + 1)} } \right)\varphi ,\nabla _A^{(n
+ 1)} \varphi } \right\rangle  = \operatorname{Re} \left\langle
{\sum\limits_{j + k = n} {\nabla _A^{(j)} \frac{{\partial A}}
{{\partial t}}} \nabla _A^{(k)} \varphi ,\nabla _A^{(n + 1)}
\varphi } \right\rangle
\\
&
 \leqslant c\left| {\nabla _M F_A } \right|\left| {\nabla _A^{(n + 1)} \varphi } \right|
 + c\left| {\nabla _M^{(n + 1)} F_A } \right|\left| {\nabla _A^{(n + 1)} \varphi } \right| + c\left| {\nabla _A^{(n + 1)} \varphi } \right|^2
 + c\left| {\nabla _A^{(n + 1)} \varphi } \right|,
\end{align*}
where ${\nabla _M F_A }$ is equal to ${\nabla _M^{(n)}
F_A }$ for the case $n=1$ and bounded for cases $n \geqslant 2$.
Thus
\begin{align}
&\frac{\partial } {{\partial t}}\left| {\nabla _A^{(n + 1)}
\varphi } \right|^2  =  - \Delta \left| {\nabla _A^{(n + 1)}
\varphi } \right|^2  - 2\left| {\nabla _A^{(n + 2)} \varphi }
\right|^2 + c\left| {\nabla _M^{(n + 1)} F_A } \right|\left|
{\nabla _A^{(n + 1)} \varphi } \right| \nonumber\\
& + c\left| {\nabla _M F_A } \right|\left| {\nabla _A^{(n + 1)}
\varphi } \right| \label{higherordervarphi} + c\left| {\nabla
_A^{(n + 1)} \varphi } \right|^2  + c\left| {\nabla _A^{(n + 1)}
\varphi } \right|.
\end{align} Similarly, from (\ref{flow2}),
\begin{align*}
&\frac{\partial }{{\partial t}}\left| {\nabla _M^{(n)} F_A }
\right|^2
 = \frac{\partial }{{\partial t}}\left| {\nabla _M^{(n)} dA} \right|^2
 = 2\left\langle {\nabla _M^{(n)} d\frac{{\partial A}}{{\partial t}},\nabla _M^{(n)} dA} \right\rangle \\
& = 2\left\langle {\nabla _M^{(n)} d\left[ { - d^* F_A  - i\operatorname{Im} \left\langle {\nabla _A \varphi ,\varphi } \right\rangle } \right],\nabla _M^{(n)} F_A }
 \right\rangle \\
& \leqslant 2\left\langle { -  \nabla _M^{(n)} \nabla_M^* \nabla_M
F_A - i\nabla _M^{(n)} d\operatorname{Im} \left\langle {\nabla _A
\varphi ,\varphi } \right\rangle ,\nabla _M^{(n)} F_A }
\right\rangle\\
&
+ c\left| {\nabla_M^{(n)}F_A } \right|^2 + c\left| {\nabla_M^{(n)}F_A } \right| \\
& \leqslant  - \Delta \left| {\nabla _M^{(n)} F_A } \right|^2  - 2\left| {\nabla _M^{(n+1)} F_A } \right|^2
- \left\langle {i\nabla _M^{(n)} d\operatorname{Im} \left\langle {\nabla _A \varphi ,\varphi } \right\rangle ,\nabla _M^{(n)} F_A } \right\rangle \\
& + c\left| {\nabla_M^{(n)}F_A } \right|^2 + c\left| {\nabla_M^{(n)}F_A } \right|, \\
\end{align*}
where we have used the Weitzenb\"ock formula
(\ref{weitzenbockforms}), the Ricci formula (\ref{Ricciformula}),
and (\ref{laplaceidentity}). Using (\ref{Fmetric}), we have
\[
\nabla _M^{(n)} d\operatorname{Im} \left\langle {\nabla _A \varphi
,\varphi } \right\rangle  = \nabla _M^{(n)} \sum\limits_{k > j}
{\left( {\left\langle {\Omega_A^{kj} \varphi ,\varphi }
\right\rangle  + 2\operatorname{Im} \left\langle {\nabla _A^j
\varphi ,\nabla _A^k \varphi } \right\rangle } \right)} dx^k
\wedge dx^j.
\]
From this and metric compatibility we find
\begin{align*}
&\left| {\left\langle {i\nabla _M^{(n)} d\operatorname{Im}
\left\langle {\nabla _A \varphi ,\varphi } \right\rangle ,\nabla
_M^{(n)} F_A } \right\rangle } \right| \leqslant \left| {i\nabla
_M^{(n)} d\operatorname{Im} \left\langle {\nabla _A \varphi
,\varphi } \right\rangle } \right|\left| {\nabla _M^{(n)} F_A }
\right|
\\
&
 \leqslant c\left| {\nabla _M^{(n)} F_A } \right|^2  + c\left| {\nabla _M^{(n)} F_A } \right|
 + c\left| {\nabla _M^{(n)} F_A } \right|\left| {\nabla _A^{(n + 1)} \varphi } \right|.
\end{align*}
Thus
 \begin{align}
\frac{\partial }{{\partial t}}\left| {\nabla _M^{(n)} F_A }
\right|^2  \leqslant & - \Delta \left| {\nabla _M^{(n)} F_A }
\right|^2  - 2\left| {\nabla _M^{(n+1)} F_A } \right|^2 + c\left|
{\nabla _M^{(n)} F_A } \right|\left| {\nabla _A^{(n + 1)} \varphi
} \right|\nonumber \\
&+c\left| {\nabla _M^{(n)} F_A } \right|^2  + c\left| {\nabla
_M^{(n)} F_A } \right|.\label{higherordereq2}
 \end{align}
Combining now equations (\ref{higherordervarphi}) and
(\ref{higherordereq2}) gives
\begin{align*}
&\quad \frac{\partial } {{\partial t}}\left( {\left| {\nabla
_A^{(n + 1)} \varphi } \right|^2  + \left| {\nabla _M^{(n)} F_A }
\right|^2
} \right) \\
&\leqslant  - \Delta \left( {\left| {\nabla _A^{(n + 1)} \varphi }
\right|^2  + \left| {\nabla _M^{(n)} F_A } \right|^2 } \right)-
2\left( {\left| {\nabla _A^{(n + 2)} \varphi } \right|^2 + \left|
{\nabla _M^{(n + 1)} F_A } \right|^2 } \right)
\\
& + c\left( {\left| {\nabla _A^{(n + 1)} \varphi } \right|^2  +
\left| {\nabla _M^{(n)} F_A } \right|^2 } \right) + c\left|
{\nabla _M^{(n + 1)} F_A } \right|\left| {\nabla _A^{(n + 1)}
\varphi } \right| \\
&+ c\left| {\nabla _M F_A } \right|\left| {\nabla _A^{(n + 1)}
\varphi } \right|+ c.
\end {align*}
Utilizing Young's inequality, we obtain (\ref{higherorder}) for
$k=n$. We now complete the proof of Lemma \ref{higherderivativeestimate}.

\underline{Case 1}. Assume $T\leq 1$. Multiplying (\ref{higherorder}) by
$e^{-c_nt}$, the maximum principle yields
\[\max_{x\in M, 0\leq t\leq T  }\left(|\nabla_A^{(n+1)} \varphi |^2 +|\nabla_M^{(n)}
 F_A|^2 \right) \leq  e^{c_n}\left ( |\nabla_A^{(n+1)} \varphi_0 |^2
 +|\nabla_M^{(n)}F_{A_0}|^2+1\right ).\]
The required result is proved.

\underline{Case 2}. Assume $T>1$.  Let $t_0$ be any time with $0\leq t_0\leq
T$. For any $t_0\leq 1$, the result follows from Case 1. For any
$t_0>1$, integrating (\ref  {higherorder}) over $M$ for $k=n-1$,
we have
\begin{align*}
 & \frac d{dt}\int_M| \nabla_A^{(n)} \varphi |^2
+|\nabla_M^{(n-1)} F_A|^2\,dV+ c'_{n-1}  \int_M |\nabla_A^{(n+1)}
\varphi |^2 +|\nabla_M^{(n)}
 F_A|^2\,dV
\\
&\leq c_{n-1}  \int_M( |\nabla_A^{(n)} \varphi |^2 +|\nabla_M^{(n-1)}
F_A|^2+1) \,dV.
\end {align*}
Integrating in $t$ on $[t_0-1,t_0]$ yields
 \begin{align*}
  \int_{t_0-1}^{t_0} \int_M |\nabla_A^{(n+1)} \varphi |^2
+|\nabla_M^{(n)}
 F_A|^2\,dV\,dt
\leq \frac{ (c_{n-1}+1)(2K^2_{n}+1)|M| }{c'_{n-1}}.
 \end {align*}
Then, using Moser's Harnack inequality in (\ref {higherorder})
with $k=n$, the required result follows. $\square$ \vspace{5mm}

\begin{Corollary} \label{local estimate} Let $(\varphi, A)$ be a solution to (\ref{flow1})-(\ref{flow2}).
Suppose $\left| {\nabla _A^{(j)} \varphi } \right| \leqslant K_n$
and $\left| {\nabla _M^{(j - 1)} F_A } \right| \leqslant K_n$ in
$P_1(x_0,t_0)$ for each $1 \leqslant j \leqslant n$ and some
constant $K_n$. Then there is a positive constant $K_{n+1}$ such
that
\[\left| {\nabla _A^{(n+1)} \varphi } \right| \leqslant
K_{n+1}, \quad \left| {\nabla _M^{(n)} F_A } \right| \leqslant
K_{n+1}  \quad\mbox{ in }  P_{1/2}(x_0,t_0).\]
\end{Corollary}
\emph{Proof.} Let $\xi$ be a smooth cut-off function
$C^{\infty}(P_1)$ satisfying $|\xi|\leq 1$ and $|\nabla
\xi|+|\Delta \xi| +|\partial_t\xi|\leq C$ in $P_{1}$ for some
constant $C>0$, and $\xi \equiv 1$ in $P_{3/4}$, $\xi \equiv 0$ on
the parabolic boundary of $P_{1}$. Multiplying  (\ref
{higherorder}) by $\xi^2$ for $k=n-1$ and integrating on $P_1$, we
have
 \begin{align*}
&c_{n-1}'\int_{P_1} \xi^2(|\nabla_A^{(n+1)} \varphi |^2
+|\nabla^{(n)}
 F_A|^2)\,dV\,dt\\
 &\leq  2\int_{P_1} (|\xi_t|+|\Delta \xi|  +|\nabla \xi|^2) (|\nabla_A^{(n)} \varphi |^2
+|\nabla^{(n-1)}
 F_A|^2)\,dV\,dt\\
 &\quad +c_{n-1}\int_{P_1}(|\nabla_A^{(n)} \varphi |^2
+|\nabla^{(n-1)}
 F_A|^2+1)  \,dV\,dt\\
 &\leq |B_1| (2K^2_n+1) (4C+2C^2+c_{n-1}).
 \end {align*}
Applying Moser's Harnack inequality to  (\ref {higherorder}) with
$k=n$ in $P_{3/4}$, the required result follows. $\square$
\vspace{5mm}

As mentioned in Section 1, we can show that concentration does not
occur in general for the Seiberg-Witten flow. We say that the energy concentrates at a point $x_0$ at time $t=T$
if there are constants $\delta$ and $R_0$ such that
\[
\mathop {\lim \sup }\limits_{t \to T} \int_{B_R (x_0 )} {\left(
{\left| {\nabla _A \varphi } \right|^2  + \left| {F_A } \right|^2
} \right)\,dV \geqslant \delta }
\]
for all $R\in (0, R_0]$. That is, as $t \to T$ we have energy
$\delta$ concentrating in smaller and smaller balls. Recall that
$\delta  > 0$ is the constant defined in Lemma
\ref{regularitytheorem}. Concentrations of amounts of energy less
than delta are ruled out by Lemma \ref{regularitytheorem}. Using Lemma
\ref{localenergyestimate} and that  the energy is bounded, it follows from the proof of
Struwe (see \cite{Harmonic} and \cite{StruweBook}) that concentration can occur at no more
than a finite number of points at $t=T$.
\begin{Lemma6} \label{blowupargument}
The energy does not concentrate at any $T\leq \infty$.
\end{Lemma6}
\emph{Proof.} We assume by contradiction that the energy
concentrates at a point $x_0$. We choose $R_0>0$ sufficiently
small so that $B_{R_0}(x_0)$ contains no concentration points
other than $x_0$. Then there exist sequences $x_m  \to x_0$, $t_m
\to T$ and a sequence of balls $B_{R_m}(x_m)$ with $R_m \to 0$
such that
\begin{align} \label{energyequalsdelta}
\delta & > \mathcal{SW}_{B_{R_m } (x_m )} (\varphi (t_m ),A(t_m
))\nonumber\\
& = \mathop {\sup }\limits_{0  \leqslant t \leqslant t_m ,\,\,x
\in B_{R_0} (x_0 )} \mathcal{SW}_{B_{R_m } (x)} (\varphi
(t),A(t))>\frac {3\delta} 4
\end{align}
for each $m$.  Choosing $C_2  =
\frac{\delta }{{4C_1 }}$,  where $C_1$ is the constant from Lemma
\ref{localenergyestimate}, and applying Lemma
\ref{localenergyestimate} to the time interval $\left[ {t,t_m }
\right]$ for some $t \in \left[ {t_m  - C_2 R_m^2 ,t_m } \right]$
gives
\begin{equation} \label{energycontradiction}
\mathcal{SW}_{B_{2R_m } (x_m )} (\varphi (t),A(t)) \geqslant
 \frac {3\delta}4 - C_1 (t_m  - t)R_m^{ - 2}  \geqslant \frac {3\delta}4  -
\frac{\delta }{4} = \frac{\delta }{2}.
\end{equation}
Define
\[
\mathcal{D}_m  = \{ (y,s):R_m y + x_m  \in B_{R_0} (x_0 ),\,s \in
\left[ { - C_2 ,0} \right]\}.
\]
Note that as $m \to \infty$, $R_m \to 0$ and $\mathcal{D}_m$ will
expand to cover $\mathbb{R}^4 \times \left[ { - C_2 ,0} \right]\
$. Similarly to the proof of Lemma \ref{regularitytheorem}, we
rescale the data to
\[
\varphi _m (y,s) = \varphi (R_m y + x_m ,R_m^2 s + t_m ),
\]
\[
A_m (y,s) = R_m A(R_m y + x_m ,R_m^2 s + t_m ),
\]
so that $\varphi _m$ and $A_m$ are defined on $\mathcal{D}_m$ and
\[
\left| {\nabla _{A_m } \varphi _m } \right|^2  = R_m^2 \left|
{\nabla _A \varphi } \right|^2,
\]
\[
\left| {F_{A_m } } \right|^2  = R_m^4 \left| {F_A } \right|^2.
\]
We next show that $R_m \varphi _m$ and $A_m$ converge locally to
  $\tilde \varphi$ and $\tilde A$ respectively, where $\tilde
\varphi$ and $\tilde A$ are defined on $\mathbb{R}^4  \times
\left[ { - C_2 ,0} \right]$. We consider the rescaled equations
\begin{equation} \label{phirescaledequation}
\frac{{\partial R_m \varphi _m }} {{\partial s}} = R_m^3
\frac{{\partial \varphi }} {{\partial t}} =  - \nabla _{A_m }^*
\nabla _{A_m } R_m \varphi _m  - \frac{1} {4}\left[ {R_m^2 S +
\left| {R_m \varphi _m } \right|^2 } \right]R_m \varphi _m,
\end{equation}
\begin{equation} \label{Arescaledequation}
\frac{{\partial A_m }} {{\partial s}} = R_m^3 \frac{{\partial A}}
{{\partial t}} =  - d^* F_{A_m }  - i\operatorname{Im}
\left\langle {\nabla _{A_m } R_m \varphi _m ,R_m \varphi _m }
\right\rangle.
\end{equation}
Note that the following argument mirrors that presented in Lemma
(\ref{regularitytheorem}) and Lemma \ref{higherderivativeestimate} for the original equations. From
(\ref{energyequalsdelta}) and Lemma \ref{regularitytheorem},
\begin{equation} \label{firstblowupbound}
R_m^2 \left| {\nabla _{A_m } \varphi _m } \right|^2  + \left|
{F_{A_m } } \right|^2  \leqslant K_1
\end{equation}
locally in $B_R(0)\times [-C_2,0]$ uniformly in $m$ where $K_1$ is
independent of $m$. Noting the similarity of these equations to
(\ref{flow1}) and (\ref{flow2}), we use (\ref{firstblowupbound})
and results identical to Lemma \ref{higherderivativeestimate} and
Corollary \ref{local estimate} to find
\begin{equation} \label{higherblowupbounds}
\left| {\nabla _{A_m }^{(n+1)} R_m \varphi _m } \right|^2  +
\left| {\nabla _M^{(n)} F_{A_m } } \right|^2  \leqslant K_{n+1}
\end{equation}
in $B_{R}(0)\times [-C_2,0]$ uniformly in $m$ for each $n \geq 0$.

If we choose our local coordinates on $B_{R_0}(x_0)$ to be orthonormal coordinates, then the metric on the rescaled space is simply $g_{ij}=\delta_{ij}$. From (\ref{timeintegralbound}), we know $\partial_tA  \in L^2( [0,\infty); L^2(M) )$ so that
\[
\int_{\mathcal{D}_m } {\left| \partial_s A_m  \right|^2 } dyds
\leqslant \int_{M \times \left[ {t_m  - C_2 R_m^2 ,t_m  } \right]}
{\left| \partial_t A  \right|^2 } dVdt \to 0.
\]
Then from (\ref {Arescaledequation}), there exists some $\tau _m  \in \left[ { - C_2,0}\right]$ such that
\[
\int_{\mathcal{D}_m (s = \tau _m )} {\left| d^* F_{A_m} \right|^2
} dy \to 0.
\]
By a result of Uhlenbeck in \cite {U2} (theorem 1.3 of \cite {U2}, see also \cite{Hong.Tian}), passing to a subsequence (without changing
notation) and in an appropriate gauge, $A_m (\tau_m)\to \tilde A$ and $R_m \varphi _m \to \tilde
\varphi$ in $C^\infty$, where   $d^*F_{\tilde A}=0$ in $\mathbb
R^4$ and $\int_{\mathbb R^4} |F_{\tilde A}|^2\,dy< C$ for some
$C>0$, and $\tilde \varphi = 0$ by the boundedness of $\varphi_m$.
Next, from (\ref{energycontradiction}),
\begin{equation} \label{middleenergy}
\int_{B_2 (0)} {R_m^2 \left| {\nabla _{A_m } \varphi _m }
\right|^2  + \left| {F_{A_m } } \right|^2  + } \frac{1} {4}R_m^4
\left| {\varphi _m } \right|^2 (S + \left| {\varphi _m } \right|^2
)dy \geqslant \frac{\delta } {2}.
\end{equation}
Since $R_m \varphi_m \to 0$, the
first and third terms of (\ref{middleenergy}) go to zero. Then we must have
\begin{equation} \label{energycontradiction2}
\int_{B_2 (0)}{\left| {F_{\tilde A } } \right|^2 } dy \geqslant \frac{\delta}{2}.
\end{equation}
We now derive a contradiction with (\ref{energycontradiction2}).
Since $F_{\tilde A}$ is harmonic in $\mathbb R^4$, the well-known
mean value formula implies that for any $x_0\in\mathbb R^4$ and
$R>0$, we have
\[ | F_{\tilde A}|(x_0)\leq \frac 1{|B_R(x_0)|}\int_{B_R(x_0)} | F_{\tilde A} |\,dy \leq \left ( \frac 1{|B_R(x_0)|}\int_{B_R(x_0)} | F_{\tilde A} |^2\,dy\right )^{1/2}\]
Letting $R\to\infty$, $ F_{\tilde A}=0$ for any $x_0\in \mathbb R^4$,
which contradicts (\ref{energycontradiction2}), as required.
 $\square$ \vspace{5mm}

Next we complete a proof of Theorem \ref{Main1}.
\vspace{5mm}

\noindent{\bf Proof of Theorem \ref{Main1}.}
By the non-concentration of the energy (Lemma
\ref{blowupargument}) at any $T \leq \infty$, there exists $R>0$ such that for any point $x \in M$ and $t \in \left[
{0,T} \right)$,
\[
\int_{B_{R} (x)} {\left( {\left| {\nabla _A \varphi } \right|^2
+ \left| {F_A } \right|^2 } \right)}(\cdot ,t)\,dV  < \delta.
\]
Then by Lemma \ref{regularitytheorem}, ${\left| {\nabla _A \varphi
} \right|^2  + \left| {F_A } \right|^2 }$ is uniformly bounded on
$P_{R/2} (x,t)$. Since $x$ and $t$ are arbitrary, ${\left|
{\nabla _A \varphi } \right|^2  + \left| {F_A } \right|^2 }$ is
uniformly bounded under the flow. From this fact and Lemma
\ref{higherderivativeestimate} we have for each $n \in \mathbb{N}$
\[
\mathop {\sup }\limits_{M \times \left[ {0,\infty } \right)}
\left( {\left| {\nabla _A^{(n)} \varphi } \right|^2  + \left|
{\nabla _M^{(n-1)} F_A } \right|^2 } \right) \leqslant K_n.
\]
Note that equations (\ref{flow1})-(\ref{flow2}) depend only on these bounded quantities. It is then elementary to show using the Sobolev embedding theorem that $(\varphi(t),A(t))$ converges to smooth data $(\varphi(T),A(T))$ as $t \to T$. In conjunction with local existence, this shows Theorem \ref{Main1}. $\square$
\vspace{5mm}

\section{Convergence}

In this section we prove Theorem 2. That is, we show that the flow (\ref{flow1})-(\ref{flow2}) converges uniquely to a critical point of the functional (\ref{swfunc}). Since convergence is only possible up to gauge, throughout this section we assume an appropriate choice of gauge. We denote a critical point of the Seiberg-Witten functional by $(\varphi_\infty, A_\infty)$, and write $\varphi=\varphi_\infty + \psi$ and $A=A_\infty +a$, where $(\varphi,A)$ denotes a solution to the flow. For simplicity, we denote $\left\| \varphi \right\|+\left\| A \right\|$ by $\left\| (\varphi,A) \right\|$ for any norm $\left\| {\, \cdot \,} \right\|$. The proof depends on the following lemmas.

\begin{Lemma9} \label{subsequence}
For each $k>0$, there exist sequences $\{ t_n \}$ and $\{g_n\} \subset \mathscr G$ with $t_n \to \infty$ such that $g_n  \cdot (\varphi (t_n ),A(t_n ))$ converges in $H^k$ to a critical point $(A_\infty, \varphi_\infty)$.
\end{Lemma9}
\begin{proof}
Integrating the energy inequality we find
\begin{equation*}
\int_0^\infty  {\left\| {\left( {\frac{{\partial \varphi }}{{\partial t}},\frac{{\partial A}}
{{\partial t}}} \right)} \right\|_{L^2 } }  \leqslant c.
\end{equation*}
It follows that there exists a sequence $\{t_n\}$ such that
\begin{equation} \label{derivativetozero}
\left\| {\left( {\frac{{\partial \varphi }}{{\partial t}}(t_n ),\frac{{\partial A}}{{\partial t}}(t_n )} \right)} \right\|_{L^2 }  \to 0.
\end{equation}
Next, recall from Lemma \ref{higherderivativeestimate} that we have uniform bounds on the quantities $\left\| {\varphi } \right\|_{H^k }$  and $\left\| {F_A } \right\|_{H^k }$ for each $k \geq 0$. It follows from a theorem of Uhlenbeck (theorem 1.3 of \cite {U2}) that in an appropriate (time varying) gauge, we also have uniform bounds on $\left\| {A} \right\|_{H^k}$ for each $n \geq 0$. For each $k \geq 0$, from the Rellich-Kondrachov theorem we can pass to a subsequence of $\{t_n\}$ (without changing notation) such that $(\varphi (t_n ),A(t_n ))$ converges in $H^k$ up to gauge to a point ${(\varphi _\infty  ,A_\infty  )}$. It remains to show that ${(\varphi _\infty  ,A_\infty  )}$ is a critical point. From (\ref{derivativetozero}) we have the required result.
\end{proof}

\begin{Lemma10} \label{continuousdependence}
On any finite time interval, the solution to the flow depends continuously on the initial conditions. That is, if $(\varphi _1 (t),A_1 (t))$ and $(\varphi _2 (t),A_2 (t))$ are two solutions to the flow with different initial values, then for any $T>0$ there exists a constant $c$ such that
\begin{equation}
\left\| {(\varphi _1 (T),A_1 (T)) - (\varphi _2 (T),A_2 (T))} \right\|_{H^k }  \leqslant c\left\| {(\varphi _1 (0),A_1 (0)) - (\varphi _2 (0),A_2 (0))} \right\|_{H^k }.
\end{equation}
\end{Lemma10}

\begin{proof}
Recall that in the gauge of theorem 1.3 of \cite{U2}, we have uniform bounds on $\varphi$, $A$, and all of their derivatives. In this gauge, we also know that $d^*A=0$. Using these facts and the expansion
\[
\nabla _A^* \nabla _A \varphi  =  - \nabla _{A_\infty  }^* \nabla _{A_\infty  } \varphi  + a\# \nabla _{A_\infty  } \varphi  + \nabla _M a\# \varphi  + a\# a\# \varphi,
\]
we can write
\begin{equation} \label{ctsdep1}
\frac{\partial }{{\partial t}}(\varphi _1  - \varphi _2 ) =  - \nabla _{A_\infty  }^* \nabla _{A_\infty  } (\varphi _1  - \varphi _2 ) + f,
\end{equation}
\begin{equation} \label{ctsdep2}
\frac{\partial }{{\partial t}}(A_1  - A_2 ) =  - \Delta (A_1  - A_2 ) + g,
\end{equation}
where $f$ and $g$ comprise the lower order terms from (\ref{flow1}) and (\ref{flow2}), and $\left\| f \right\|_{H^k }$ and $\left\| g \right\|_{H^k }$ are both bounded by $c\left\| {(\varphi _1  - \varphi _2 ,A_1  - A_2) } \right\|_{H^k }$. When $f=g=0$, we simply have the heat equation, whose solution depends continuously on the intial data in the $H^k$ norm. When the data is small in the $H^k$ norm, the perturbations $f$ and $g$ will be small in the $H^k$ norm also. Thus $\varphi$ and $A$ depend continuously on their initial values.
\end{proof}

\begin{Lemma11} \label{Lojasiewicz}
Let $(\varphi_\infty,A_\infty)$ be a critical point of the Seiberg-Witten functional. There exist constants $\varepsilon_1>0$, $\frac{1}{2}<\gamma<1$ and $c>0$ such that if
\[
\left\|(\varphi,A)-(\varphi_\infty,A_\infty)\right\|_{H^1} \leq \varepsilon_1,
\]
then
\begin{equation}
\left\| {\left( {\frac{{\partial \varphi }}{{\partial t}},\frac{{\partial A}}{{\partial t}}} \right)} \right\|_{L^2 }  \geq c\left| {\mathcal{SW}(\varphi ,A) - \mathcal{SW}(\varphi _\infty  ,A_\infty  )} \right|^\gamma.
\end{equation}
\end{Lemma11}

\begin{proof}
The proof of this lemma is analogous to that of proposition 7.2 in \cite{Rade} and proposition 3.5 \cite{Wilkin}. While Wilkin considers in section 3 of \cite{Wilkin} the Yang-Mills-Higgs functional, he allows in the proof of this
lemma a very general functional $f: Q \to \mathbb{R}$, where $Q$
is a Hilbert manifold and $f$ is invariant under the action of
some gauge group $\mathscr G$. To apply Proposition 3.20 of
\cite{Wilkin}, we need only check that the operator
\[
H_{\mathcal{SW}} + \rho  _\infty \rho  _\infty ^* :T_\infty H \to T_\infty H
\]
is elliptic. Here $H_{\mathcal{SW}}$ represents the Hessian of the
Seiberg-Witten functional at the point
$(\varphi_\infty,A_\infty)$, and $\rho _\infty:\mathrm{Lie}(\mathscr G) \to
T_\infty M$ is the infinitesimal action of the gauge group
$\mathscr G$. The operator $ \rho _\infty ^*$ is defined by
\[
\left\langle {\rho _\infty ^* X,u} \right\rangle _{\mathrm{Lie}(\mathscr G)}  = \int_M {\left\langle {X, \rho _\infty  u} \right\rangle },
\]
for $X \in T_\infty H$ and $u \in \mathrm{Lie}(\mathscr G)$. We begin by computing the operator
\[
H_{\mathcal{SW}} (\psi ,a) = \left. {\frac{d}{{ds}}} \right|_{s = 0} Grad (\mathcal{SW}) (s\psi ,sa)
\]
where $\mathrm{Grad}(\mathcal{SW})$ represents the gradient operator of the Seiberg-Witten functional. There holds
\[
\left. {\frac{\partial }
{{\partial s}}} \right|_{s = 0} \left( {\nabla _{A_\infty   + sa}^* \nabla _{A_\infty   + sa} (\varphi _\infty   + s\psi ) + \frac{S}
{4}(\varphi _\infty   + s\psi ) + \frac{1}
{4}\left| {\varphi _\infty   + s\psi } \right|^2 (\varphi _\infty   + s\psi )} \right)
\]
\[
 = \nabla _{A_\infty  }^* \nabla _{A_\infty  } \psi  + \frac{S}
{4}\psi  + \frac{1}
{4}\left| {\varphi _\infty  } \right|^2 \psi  + \frac{1}
{2}\operatorname{Re} \left\langle {\varphi _\infty  ,\psi } \right\rangle \varphi _\infty
\]
\[
 + \left\langle {\frac{1}
{2}i\operatorname{Im} \left\langle {\varphi _\infty  ,\nabla _{A_\infty  } \psi } \right\rangle  + \frac{1}
{2}i\operatorname{Im} \left\langle {\psi ,\nabla _{A_\infty  } \varphi _\infty  } \right\rangle ,a} \right\rangle,
\]
where we use a relationship analogous to that in equation (\ref{imaginary}). Similarly,
\[
\left. {\frac{\partial }
{{\partial s}}} \right|_{s = 0} \frac{1}{2}d^* d(A_\infty   + sa) + \frac{1}{2}i\operatorname{Im} \left\langle {\nabla _{A_\infty  } (\varphi _\infty   + s\psi ) + \frac{1}
{2}sa(\varphi _\infty   + s\psi ),(\varphi _\infty   + s\psi )} \right\rangle
\]
\[
 = \frac{1}{2}d^* da + \frac{1}{2}i\operatorname{Im} \left\langle {\nabla _{A_\infty  } \varphi _\infty  ,\psi } \right\rangle  + \frac{1}{2}i\operatorname{Im} \left\langle {\nabla _{A_\infty  } \psi ,\varphi _\infty  } \right\rangle  + \frac{1}
{2}a\left| {\varphi _\infty  } \right|^2.
\]
Using the above and recalling the calculations in the proof of
Lemma \ref{energyinequalitylemma}, the Hessian at the point
$(\varphi_ \infty,A_ \infty)$ is given by
\begin {align} \label{Hessian}
\nonumber H_{\mathcal{SW}} (\psi ,a) = &\left( \nabla _{A_\infty
}^* \nabla _{A_\infty  } \psi  + \frac{1} {2}\operatorname{Re}
\left\langle {\varphi _\infty  ,\psi } \right\rangle \varphi
_\infty  + \frac{S} {4}\psi  + \frac{1}{4}\left| {\varphi _\infty
} \right|^2 \psi \right.,
\\
&\left. \frac{1}{2}d^* da + \frac{1}{4}\left| {\varphi _\infty  }
\right|^2 a  + i\operatorname{Im} \left\langle {\varphi _\infty  ,
\nabla _{A_\infty  } \psi } \right\rangle  + i\operatorname{Im}
\left\langle {\psi ,\nabla _{A_\infty  } \varphi _\infty  }
\right\rangle \right). \end {align}

In the following, we continue to use $(\cdot,\cdot)$ to denote an
element of the configuration space, i.e., $(\psi ,a) \in \Gamma (S^
+  ) \times \mathscr A$. Now, note that if $g(t)$ represents a
path through the gauge group $\mathscr G$ with $g(0)=I$, then
\begin{align*}
& \rho _\infty  (g'(0)) = \frac{1}{{\sqrt 2 }} \left. {\frac{\partial }{{\partial t}}} \right|_{t = 0} (g(t)^* (\varphi _\infty  ,A_\infty  )) \\
&  = \frac{1}{{\sqrt 2 }} \left. {\frac{\partial }{{\partial t}}} \right|_{t = 0} (g(t)^{-1}\varphi _\infty  ,A_\infty   + 2g(t)^{ - 1} dg(t))\\
&  = \frac{1}{{\sqrt 2 }} ( - g'(0)\varphi _\infty  ,2dg'(0)),
\end{align*}
where we write $g(t)=e^{i \theta}$ for some function $\theta
(t,x)$ defined locally on the manifold $M$, so that
\[
\left. {\frac{\partial } {{\partial t}}} \right|_{t = 0} 2g(t)^{ -
1} dg(t) = \left. {\frac{\partial } {{\partial t}}} \right|_{t =
0} 2id\theta  = 2dg'(0).
\]
It follows that
\begin{align*}
& \left\langle { \rho _\infty ^* (\psi ,a),g'(0)} \right\rangle _{\mathrm{Lie}(\mathscr G)}  = \frac{1}{{\sqrt 2 }} \int_M {\left\langle {\psi , - g'(0)\varphi _\infty  } \right\rangle  + \left\langle {a,2dg'(0)} \right\rangle } \\
& = \frac{1}{{\sqrt 2 }}\left\langle {\left\langle {\psi ,\varphi
_\infty  } \right\rangle ,g'(0)} \right\rangle _{\mathrm{Lie}(\mathscr G)}
+ \frac{1}{{\sqrt 2 }} \left\langle {\frac{1} {2}d^* a,g'(0)}
\right\rangle _{\mathrm{Lie}(\mathscr G)},
\end{align*}
that is,
\[
\rho  _\infty ^* (\psi ,a) = \frac{1}{{\sqrt 2 }}\left(
{\left\langle {\psi ,\varphi _\infty  } \right\rangle ,\frac{1}
{2}d^* a} \right)
\]
and
 \begin {equation} \label{rhorhostar}
\rho  _\infty \rho _\infty ^* (\psi ,a) = \frac{1}{2} \left( { -
\left\langle {\psi ,\varphi _\infty  } \right\rangle \varphi
_\infty  ,dd^* a} \right).
\end{equation}
Comparing (\ref{rhorhostar}) with (\ref{Hessian}), we find that
$H_{\mathcal{SW}}  + \rho _\infty \rho _\infty ^*$ is an elliptic
operator, as required. Then, the required result follows from the same arguments as for Theorem 3.19 in \cite{Wilkin}.
\end{proof}

\begin{Lemma12} \label{interior}
There exists a constant $c$ such that if $T \geq 0$ and $S>1$ are such that $0 \leq T \leq S-1$, then
\begin{equation}
\int_{T + 1}^S {\left\| {\left( {\frac{{\partial \varphi }}
{{\partial t}},\frac{{\partial A}}
{{\partial t}}} \right)} \right\|_{H^k } }  \leqslant c\int_T^S {\left\| {\left( {\frac{{\partial \varphi }}
{{\partial t}},\frac{{\partial A}}
{{\partial t}}} \right)} \right\|_{L^2 } }.
\end{equation}
\end{Lemma12}

\begin{proof}
We define $G = (G_1 ,G_2 ) = \left( {\frac{{\partial \varphi }}{{\partial t}},\frac{{\partial A}}{{\partial t}}} \right)$. Noting that
\[
\nabla _A^* \nabla _A \varphi  =  - \nabla _{A_\infty  }^* \nabla _{A_\infty  } \varphi  + a\# \nabla _{A_\infty  } \varphi  + \nabla _M a\# \varphi  + a\# a\# \varphi,
\]
we have
\[
\frac{{\partial G_1 }}
{{\partial t}} =  - \nabla _{A_\infty  }^* \nabla _{A_\infty  } G_1  + G_2 \# \nabla _{A_\infty  } \varphi  + a\# \nabla _{A_\infty  } G_1  + \nabla _M G_2 \# \varphi  + \nabla _M a\# G_1
\]
\[
 + G_2 \# a\# \varphi  + a\# a\# G_1  - \frac{S}{4}G_1  + \varphi \# \varphi \# G_1,
\]
and
\[
\frac{{\partial G_2 }}{{\partial t}} =  - d^* dG_2  + G_2 \# \varphi \# \varphi  + a\# \varphi \# G_1  + \nabla _{A_\infty  } \varphi \# G_1  + \varphi \# \nabla _{A_\infty  } G_1.
\]
Using the Bianchi identify and the Weitzenb\"ock formula
(\ref{weitzenbockforms}) we can write
\begin{align*}
 - d^* dG_2 & =  - \Delta G_2  - idd^*\operatorname{Im} \left\langle {\nabla _A \varphi ,\varphi } \right\rangle \\
&  =  - \nabla _M^* \nabla _M G_2  - idd^*\operatorname{Im} \left\langle {\nabla _A \varphi ,\varphi } \right\rangle  + R_M \# G_2.
\end{align*}
Using metric compatibility and equation (\ref{flow1}), we compute in normal coordinates
\begin{align*}
 &- idd^* \operatorname{Im} \left\langle {\nabla_A \varphi ,\varphi } \right\rangle
 = id*d\left( {\operatorname{Im} \left\langle {\nabla _A^j \varphi ,\varphi } \right\rangle *dx^j } \right) \\
&  = id*\left( {\left[ {\operatorname{Im} \left\langle {\nabla _A^k \nabla _A^j \varphi ,\varphi } \right\rangle
+ \operatorname{Im} \left\langle {\nabla _A^j \varphi ,\nabla _A^k \varphi } \right\rangle } \right]dx^k  \wedge *dx^j } \right)\\
&  =  id*\left( {\operatorname{Im} \left\langle {\nabla _A^j \nabla _A^j \varphi ,\varphi } \right\rangle dV} \right)
 = -id\operatorname{Im} \left\langle {\nabla _A^* \nabla _A \varphi ,\varphi } \right\rangle \\
&  = -i\operatorname{Im} \left( {\left\langle {\nabla_A \nabla_A^* \nabla _A \varphi ,\varphi } \right\rangle  + \left\langle {\nabla_A^* \nabla_A \varphi ,\nabla_A \varphi } \right\rangle } \right) \\
&  = -i\frac{1}{4}\left[ {S + \left| \varphi  \right|^2 } \right]\operatorname{Im} \left( {\left\langle {\nabla _A \varphi ,\varphi } \right\rangle  + \left\langle {\varphi ,\nabla _A \varphi } \right\rangle } \right) \\
& \quad- i\frac{1}{4}\operatorname{Im} \left\langle {\left[ {dS + d\left| \varphi  \right|^2 } \right]\varphi ,\varphi } \right\rangle \\
& \quad + \varphi \# \nabla _{A_\infty  } G_1  + \nabla _{A_\infty  } \varphi \# G_1  + a\# \varphi \# G_1 \\
&  = \varphi \# \nabla _{A_\infty  } G_1  + \nabla _{A_\infty  } \varphi \# G_1  + a\# \varphi \# G_1,
\end{align*}
where the second term in line two and the first two terms in the second to last expression are zero.
Thus recalling the uniform bounds on $\varphi$, $A$ and their derivatives (see the proof of Lemma \ref{subsequence}), we can combine all of the above in the compact form
\[
\frac{{\partial G}}{{\partial t}} + \nabla ^* \nabla G = V_0 \# G + V_1 \# \nabla G,
\]
where the $V_j$ are smooth vectors having all derivatives uniformly bounded, and $\nabla$ acts as $\nabla_{A_\infty}$ on sections of $\mathcal S^+$ and as $\nabla_M$ on forms. This equation is of the same form as the equation in the proof of Proposition 3.6 of \cite{Wilkin}, and the rest of the proof is the same. Note that since we have uniform bounds on all derivatives, we do not need to require the assumption $\left\| (\varphi(T),A(T))-(\varphi_\infty,A_\infty) \right\|_{H^k}<\varepsilon$ as in \cite{Rade} and \cite{Wilkin}.
\end{proof}

\begin{Lemma13} \label{uniqueness}
Let $k>0$. There exists $\varepsilon>0$ such that if for some $T>0$
\begin{equation}\label{48}
\left\| (\varphi(T),A(T))-(\varphi_\infty,A_\infty)
\right\|_{H^k}<\varepsilon,
\end{equation}
then either $\mathcal{SW}(\varphi(t),A(t))<\mathcal{SW}(\varphi_\infty,A_\infty)$
for some $t>T$, or $(\varphi(t),A(t))$ converges in $H^k$ to a critical point $(\varphi_\infty',A_\infty')$
where $\mathcal{SW}(\varphi_\infty',A_\infty')=\mathcal{SW}(\varphi_\infty,A_\infty)$ and
\begin{equation}\label{estimate3}
\left\| {(\varphi _\infty '  ,A_\infty ' ) - (\varphi _\infty
,A_\infty  )} \right\|_{H^k }   \leqslant c\left\| {(\varphi
(T),A(T)) - (\varphi _\infty  ,A_\infty  )} \right\|_{H^k }^{2(1 -
\gamma )}
\end{equation}
where $\gamma$ is as in Lemma \ref{Lojasiewicz}. We also have the following
convergence estimate:
\begin{equation} \label{estimate4}
\left\| {(\varphi (t),A(t)) - (\varphi _\infty ' ,A_\infty ' )} \right\|_{H^k }  \leqslant c(t-T)^{ - (1 - \gamma )/(2\gamma  - 1)}.
\end{equation}

\end{Lemma13}

\begin{proof}
We set
\begin{displaymath}
\Delta \mathcal{SW}(t) = \mathcal{SW}(\varphi (t),A(t)) -
\mathcal{SW}(\varphi _\infty  ,A_\infty).
\end{displaymath}
Then, we can assume that $\Delta \mathcal{SW}(t) \geq 0$ for all
$t$. Otherwise, the required result is proved.

We note
\begin{align*}
 &\quad \int_M {\left| {\nabla _A \varphi }
\right|^2  - \left| {\nabla _{A_\infty  } \varphi _\infty  }
\right|} ^2
 \\
 &
= \int_M {\left| {\nabla _{A_\infty  } \psi  + \frac{1}
{2}a\varphi _\infty   + \frac{1} {2}a\psi } \right|^2 }  +
2\operatorname{Re} \left\langle {\nabla _{A_\infty  } \varphi
_\infty  ,\nabla _{A_\infty  } \psi  + \frac{1} {2}a\varphi
_\infty   + \frac{1} {2}a\psi } \right\rangle
 \\
 &
= \int_M {\left| {\nabla _{A_\infty  } \psi } \right|^2  + \left|
{\frac{1} {2}a\varphi _\infty  } \right|^2  + \left| {\frac{1}
{2}a\psi } \right|^2 }  + 2\operatorname{Re} \left\langle {\nabla
_{A_\infty  } \psi ,\frac{1} {2}a\varphi _\infty  } \right\rangle
 \\
 &
\quad + 2\operatorname{Re} \left\langle {\nabla _{A_\infty  } \psi
,\frac{1} {2}a\psi } \right\rangle  + 2\operatorname{Re}
\left\langle {\frac{1} {2}a\varphi _\infty  ,\frac{1} {2}a\psi }
\right\rangle + 2\operatorname{Re} \left\langle {\nabla _{A_\infty
} \varphi _\infty  ,\nabla _{A_\infty  } \psi } \right\rangle
 \\
 &
 \quad  +
2\operatorname{Re} \left\langle {\nabla _{A_\infty  } \varphi
_\infty  ,\frac{1} {2}a\varphi _\infty  } \right\rangle  +
2\operatorname{Re} \left\langle {\nabla _{A_\infty  } \varphi
_\infty  ,\frac{1} {2}a\psi } \right\rangle
 \\
 &
= \int_M {\left| {\nabla _{A_\infty  } \psi } \right|^2  +
\frac{1} {4}\left| {\varphi _\infty  } \right|^2 \left| a
\right|^2  + \frac{1} {4}\left| \psi  \right|^2 \left| a \right|^2
}  + \left\langle {a,i\operatorname{Im} \left\langle {\nabla
_{A_\infty  } \psi, \varphi _\infty} \right\rangle } \right\rangle
 \\
 &\quad + \left\langle {a,i\operatorname{Im} \left\langle {\nabla
_{A_\infty  } \psi, \psi} \right\rangle } \right\rangle  +
\frac{1} {2}\left| a \right|^2 \operatorname{Re} \left\langle
{\varphi _\infty  ,\psi } \right\rangle + 2\operatorname{Re}
\left\langle {\nabla _{A_\infty  }^* \nabla _{A_\infty  } \varphi
_\infty  ,\psi } \right\rangle
\\
 &\quad + \left\langle {a,i\operatorname{Im} \left\langle {\nabla _{A_\infty  } \varphi _\infty, \varphi _\infty } \right\rangle }
  \right\rangle  + \left\langle {a,i\operatorname{Im} \left\langle {\nabla _{A_\infty} \varphi _\infty , \psi } \right\rangle } \right\rangle,
\end{align*}
 where we again use a relationship analogous to that in
equation (\ref{imaginary}). It is easy to see
\[
\int_M {\frac{1}
{2}\left| {F_A } \right|^2  - \frac{1}
{2}\left| {F_{A_\infty  } } \right|^2 }  = \int_M {\frac{1}
{2}\left| {da} \right|^2 }  + \left\langle {dA_\infty, da} \right\rangle.
\]
We have also
\[
\int_M {\frac{S} {4}\left| \varphi  \right|^2  - \frac{S}
{4}\left| {\varphi _\infty  } \right|^2 }  = \int_M {\frac{S}
{4}\left( {\left| \psi  \right|^2  + 2\operatorname{Re}
\left\langle {\psi ,\varphi _\infty  } \right\rangle } \right)}
\]
and
 \begin{align*}
&\frac{1} {8}\int_M {\left| \varphi  \right|^4  - \left| {\varphi
_\infty  } \right|^4 }
 = \frac{1}{8}\int_M {\left( {\left| {\varphi _\infty  } \right|^2  + \left| \psi  \right|^2
  + 2\operatorname{Re} \left\langle {\varphi _\infty  ,\psi } \right\rangle } \right)^2  - \left| {\varphi _\infty  } \right|^4 }
\\&
 = \frac{1}
{8}\int_M { {\left| \psi  \right|^4  + 4\operatorname{Re} \left\langle {\varphi _\infty  ,\psi } \right\rangle ^2  + 2\left| {\varphi _\infty  }
 \right|^2 \left| \psi  \right|^2  + 4\left| {\varphi _\infty  } \right|^2 \operatorname{Re} \left\langle {\varphi _\infty  ,\psi } \right\rangle } }
\\
&\quad    +4\left| \psi  \right|^2 \operatorname{Re} \left\langle
{\varphi _\infty  ,\psi } \right\rangle.
 \end{align*}
Recalling that $(\varphi_\infty,A_\infty)$ satisfies the critical
point equations (\ref{eq1}) and (\ref{eq2}), we have
\[
\int_M {\left\langle {\nabla _{A_\infty  }^* \nabla _{A_\infty  } \varphi _\infty   + \frac{S}
{4}\varphi _\infty   + \frac{1}
{4}\left| {\varphi _\infty  } \right|^2 \varphi _\infty  ,\psi } \right\rangle }  = 0
\]
and
\[
\int_M {\left\langle {d^* dA_\infty   + i\operatorname{Im} \left\langle {\nabla _{A_\infty  } \varphi _\infty  ,\varphi _\infty  } \right\rangle ,a} \right\rangle }  = 0,
\]

Combining above estimates, we have
\[
\Delta \mathcal{SW}(t) = \int_M {\left| {\nabla _{A_\infty  } \psi
} \right|^2  + \frac{1} {2}\operatorname{Re} \left\langle {\varphi
_\infty  ,\psi } \right\rangle ^2  + \frac{S} {4}\left| \psi
\right|^2 }  + \frac{1}{4}\left| {\varphi _\infty  } \right|^2
\left( {\left| \psi  \right|^2  + \left| a \right|^2 } \right)
\]
\begin{equation} \label{functionaldiff}
 + \frac{1}
{2}\left| {da} \right|^2   + \left\langle {a,i\operatorname{Im}
\left\langle {\nabla _{A_\infty  } \psi, \varphi _\infty}
\right\rangle } \right\rangle  + \left\langle
{a,i\operatorname{Im} \left\langle {\nabla _{A_\infty  } \varphi
_\infty, \psi } \right\rangle } \right\rangle  + O(3),
\end{equation}
where
\[
O(3) = \int_M {\frac{1} {2}\left( {\left| \psi  \right|^2  +
\left| a \right|^2 } \right)\operatorname{Re}  \left\langle
{\varphi _\infty  ,\psi } \right\rangle  + \left\langle
{a,i\operatorname{Im} \left\langle {\nabla _{A_\infty  } \psi,
\psi} \right\rangle } \right\rangle }  + \frac{1} {8}\left| \psi
\right|^4  + \frac{1} {4}\left| \psi  \right|^2 \left| a
\right|^2.
\]

Since $\Delta \mathcal{SW}(t)$ is a polynomial functional and
$(A_\infty, \varphi_\infty)$ is a critical point, the first order
terms of $\Delta \mathcal{SW}(t)$ vanish. Then for $\varepsilon$ small enough we have
\begin{align} \label{differencebound}
\Delta \mathcal{SW}(T) & \leqslant c\left\| {(\varphi (T),A(T)) - (\varphi _\infty  ,A_\infty  )} \right\|_{H^1 }^2 \nonumber \\
& \leqslant c\left\| {(\varphi (T),A(T)) - (\varphi _\infty  ,A_\infty  )} \right\|_{H^k }^2.
\end{align}
From the continuous dependence on initial conditions (Lemma
\ref{continuousdependence}), for $\varepsilon$ in (\ref{48})
sufficiently small we have for $t \in [T,T+1]$,
\begin{equation} \label{ctsdepend}
\left\| {(\varphi (t),A(t)) - (\varphi _\infty  ,A_\infty  )} \right\|_{H^k }  < \frac{1}
{2}\varepsilon _1.
\end{equation}
We claim that if $\varepsilon$ is sufficiently small, then for all $t \geq T$ we have
\begin{equation}
\left\| {(\varphi (t),A(t)) - (\varphi _\infty  ,A_\infty  )} \right\|_{H^k }  < \varepsilon _1.
\end{equation}
Suppose by contradiction that $S>T$ is the smallest number such that \\ $\left\| {(\varphi (S),A(S)) - (\varphi _\infty  ,A_\infty  )} \right\|_{H^k }  \geqslant \varepsilon _1$. From Lemma (\ref{Lojasiewicz}) we have for $T \leq t \leq S$,
\begin{align} \label{integrate}
\frac{d}{{dt}}(\Delta \mathcal{SW}(t))^{1 - \gamma }
&= - c(1 - \gamma )(\Delta \mathcal{SW}(t))^{- \gamma } \left\| {\left( {\frac{{\partial \varphi }}{{\partial t}},\frac{{\partial A}}{{\partial t}}} \right)} \right\|_{L^2 }^2  \nonumber\\
& \leq  - c\left\| {\left( {\frac{{\partial \varphi }}{{\partial t}},\frac{{\partial A}}
{{\partial t}}} \right)} \right\|_{L^2 }.
\end{align}
Integrating (\ref{integrate}) in time gives
\begin{equation} \label{afterintegration}
\int_T^S {\left\| {\left( {\frac{{\partial \varphi }}{{\partial t}},\frac{{\partial A}}
{{\partial t}}} \right)} \right\|_{L^2 } } \leq c\Delta (\mathcal{SW}(T))^{1 - \gamma }.
\end{equation}
Recalling (\ref{differencebound}),
\begin{equation} \label{integrated}
\int_T^S {\left\| {\left( {\frac{{\partial \varphi }}
{{\partial t}},\frac{{\partial A}}
{{\partial t}}} \right)} \right\|_{L^2 } }  \leqslant c\left\| {(\varphi (T),A(T)) - (\varphi _\infty  ,A_\infty  )} \right\|_{H^k }^{2(1 - \gamma )}  \leqslant c\varepsilon ^{2(1 - \gamma )}.
\end{equation}
From (\ref{ctsdepend}) we know that $S>T+1$, and then
\[
\int_{T + 1}^S {\left\| {\left( {\frac{{\partial \varphi }}
{{\partial t}},\frac{{\partial A}}
{{\partial t}}} \right)} \right\|_{H^k }  \geqslant } \left\| {\int_{T + 1}^S {\left( {\frac{{\partial \varphi }}
{{\partial t}},\frac{{\partial A}}
{{\partial t}}} \right)} } \right\|_{H^k }
\]
\[
 \geqslant \left\| {(\varphi (S),A(S)) - (\varphi _\infty  ,A_\infty  )} \right\|_{H^k }  - \left\| {(\varphi (T + 1),A(T + 1)) - (\varphi _\infty  ,A_\infty  )} \right\|_{H^k }
\]
\[
 \geqslant \varepsilon _1  - \frac{1}{2}\varepsilon _1.
\]
Then using our results above and Lemma \ref{interior}, we find
\[
\frac{1}{2}\varepsilon _1  \leqslant c\varepsilon ^{2(1 - \gamma )},
\]
which is impossible for $\varepsilon$ small enough. Thus, as claimed, for $\varepsilon$ small enough we have $\left\| {(\varphi (t),A(t)) - (\varphi _\infty  ,A_\infty  )} \right\|_{H^k }  < \varepsilon _1$ for all $t \geq T$.

Finally, letting $S \to \infty$ in Lemma \ref{interior} and (\ref{afterintegration}) we have
\begin{align} \label{showconvergence}
\int_{t_1 + 1}^\infty  {\left\| {\left( {\frac{{\partial \varphi }}{{\partial t}},\frac{{\partial A}}
{{\partial t}}} \right)} \right\|_{H^k } }
& \leq c\int_{t_1}^\infty  {\left\| {\left( {\frac{{\partial \varphi }}{{\partial t}},\frac{{\partial A}}
{{\partial t}}} \right)} \right\|_{L^2}}  \nonumber\\
& \leq  c(\Delta \mathcal{SW}(t_1))^{1 - \gamma }
\end{align}
for any $t_1 \geq T$. From Lemma \ref{subsequence} and Lemma \ref{Lojasiewicz}, we have
\[
\int_{t_1}^\infty  {\left\| {\left( {\frac{{\partial \varphi }}{{\partial t}},\frac{{\partial A}}
{{\partial t}}} \right)} \right\|_{H^k } } \to 0
\]
as $t_1 \to \infty$. This establishes unique convergence of the flow in the $H^k$ norm to a point $(\varphi' _\infty  ,A'_\infty  )$,
provided that $\left\| (\varphi(T),A(T))-(\varphi_\infty,A_\infty) \right\|_{H^k}<\varepsilon$ for some $T$.

As in Lemma \ref{continuousdependence}, it follows that $(\varphi'
_\infty  ,A'_\infty)$ is a critical point, and it follows from Lemma \ref{Lojasiewicz} that
$\mathcal{SW}(\varphi_\infty',A_\infty')=\mathcal{SW}(\varphi_\infty,A_\infty)$. Then from (\ref{showconvergence}) and (\ref{differencebound}) we
have
\begin{align} \label{anestimate}
&\left\| {(\varphi (T + 1),A(T+ 1))-(\varphi' _\infty  ,A'_\infty  )} \right\|_{H^k } \nonumber \\
& \leq \int_{T + 1}^\infty  {\left\| {\left( {\frac{{\partial \varphi }}{{\partial t}},\frac{{\partial A}}
{{\partial t}}} \right)} \right\|_{H^k } }  \leq c(\Delta \mathcal{SW}(T))^{1 - \gamma }  \\
& \leq c\left\| {(\varphi (T),A(T)) - (\varphi _\infty  ,A_\infty  )} \right\|_{H^k }^{2(1 - \gamma )} \nonumber\\
& \leq c\left\| {(\varphi (T),A(T)) - (\varphi _\infty  ,A_\infty  )} \right\|_{H^k } \nonumber
\end{align}
since $\gamma \in (\frac{1}{2},1)$, and from Lemma \ref{continuousdependence},
\[
\left\| {(\varphi (T + 1),A(T + 1)) - (\varphi _\infty  ,A_\infty  )} \right\|_{H^k }  \leqslant c\left\| {(\varphi (T),A(T)) - (\varphi _\infty  ,A_\infty  )} \right\|_{H^k}.
\]
The estimate (\ref{estimate3}) follows from the above two inequalities. It remains to show (\ref{estimate4}). As in (\ref{anestimate}), for $t \geq T$ we have
\[
\left\| {(\varphi (t + 1),A(t + 1))-(\varphi' _\infty  ,A'_\infty  )} \right\|_{H^k }  \leqslant c(\Delta \mathcal{SW}(t))^{1 - \gamma }.
\]
Then from Lemma \ref{Lojasiewicz} we have
\[
\frac{d}
{{dt}}\Delta \mathcal{SW}(t) =  - c\left\| {\left( {\frac{{\partial \varphi }}
{{\partial t}},\frac{{\partial A}}
{{\partial t}}} \right)} \right\|_{L^2 }^2  \leqslant  - c(\Delta SW(t))^{2\gamma },
\]
which implies that
\begin{equation} \label{functionalestimate}
\Delta \mathcal{SW}(t) \leqslant c(t - T)^{ - 1/(2\gamma  - 1)}.
\end{equation}
Thus combining the above, for $t \geq T+1$ we find
\begin{equation} \label{dropconstant}
\left\| {(\varphi (t),A(t))-(\varphi' _\infty  ,A'_\infty  )} \right\|_{H^k }  \leqslant c(t-T-1)^{-(1 - \gamma)/(2 \gamma -1) }.
\end{equation}
Note that since the left-hand side is bounded under the flow, by adjusting the constant $c$ if necessary, we can drop the constant $1$, and (\ref{estimate4}) follows.
\end{proof}

We now complete the proof of Theorem \ref{Main2}.

\begin{proof}[Proof of Theorem \ref{Main2}]
From the convergence of a subsequence $\{t_k\}$ of the flow to a
critical point $(\varphi_\infty ,A_\infty )$ (Lemma
\ref{subsequence}),  we know that there exists a $T$ such that
$\left\| (\varphi(T),A(T))-(\varphi_\infty,A_\infty)
\right\|_{H^k}<\varepsilon$. We can then apply Lemma \ref{uniqueness}. Note that in deriving (\ref{finalestimate}), as for (\ref{dropconstant}), by adjusting the constant $c$ if necessary we can drop the constant $T$.

Finally, we show that the limit depends continuously on
the initial data in the space $\{(\varphi_0, A_0): SW (\varphi
(t),A(t)) \to \lambda \}$ as $t\to \infty$.  Let ${(\varphi
(t),A(t))}$ be a solution to the flow which converges to
${(\varphi _\infty ,A_\infty  )}$ as $t \to \infty$. Let
${(\varphi '(t),A'(t))}$ be another solution to the flow with
initial data $(\varphi' (0),A'(0))$ with
\[\mathop {\lim
}\limits_{t \to \infty } \mathcal{SW}(\varphi '(t),A'(t))
=\mathcal{SW}(\varphi' _\infty ,A'_\infty  )= \mathcal{SW}(\varphi
_\infty ,A_\infty  ).\]
From Lemma
\ref{uniqueness}, for any $\beta_1 > 0$ there exists a $\beta_2 >
0$ such that if for some $T \geq 0$,
\[
\left\| {(\varphi' (T),A'(T)) - (\varphi _\infty  ,A_\infty  )}
\right\|_{H^k }  \leqslant \beta _2,
\]
then
$(\varphi' (t),A'(t))$ converges in $H^k$ as $t \to \infty$ to a
critical point ${(\varphi' _\infty  ,A'_\infty  )}$, and further $\left\|
{(\varphi' _\infty ,A'_\infty ) - (\varphi _\infty ,A_\infty )}
\right\|_{H^k } \leqslant \beta _1$. Choose $T$ such that
\[
\left\| {(\varphi (T),A(T)) - (\varphi _\infty  ,A_\infty  )}
\right\|_{H^k }  \leqslant \frac{{\beta _2 }}{2}.
\]
From Lemma {\ref{continuousdependence}}, there exists $\beta_3>0$
such that if
\[\left\| {(\varphi' (0),A'(0)) - (\varphi (0),A(0))}
\right\|_{H^k }  \leqslant \beta_3,\] then $\left\| {(\varphi
(T),A(T)) - (\varphi' (T),A'(T))} \right\|_{H^k }  \leqslant
\frac{{\beta _2 }}{2}$. Applying the triangle inequality, for any
$\beta_1>0$ there exists a $\beta_3>0$ such that if
\[\left\|
{(\varphi' (0),A'(0)) - (\varphi (0),A(0))} \right\|_{H^k } \leq
\beta_3,\]
then \[\left\| {(\varphi' _\infty  ,A'_\infty  ) -
(\varphi _\infty  ,A_\infty  )} \right\|_{H^k }   \leq \beta _1.\]
This completes the proof of Theorem \ref{Main2}.
\end{proof}

\section{Perturbed Functional}

One can also consider the perturbed Seiberg-Witten equations
\begin{equation} \label{sweqpert}
D_A \varphi  = 0,\;\;\;\;\;\;\;\;\;\;\;\;\;F_A^ +   = \frac{1}
{4}\left\langle {e_j e_k \varphi ,\varphi } \right\rangle e^j
\wedge e^k + \mu
\end{equation}
and the corresponding perturbed Seiberg-Witten functional
\[
\mathcal{SW}_{\mu}(\varphi ,A) = \int_M {\left| {D_A \varphi }
\right|^2 }  + \left| {F_A^ +   - \frac{1} {4}\left\langle {e_j
e_k \varphi ,\varphi } \right\rangle e^j  \wedge e^k - \mu }
\right|^2.
\]
\begin{equation} \label{swfuncpert}
= \int_M {\left| {\nabla _A \varphi } \right|^2  + \left| {F_A^ +
} \right|^2  + \frac{S} {4}\left| \varphi  \right|^2  + \frac{1}
{8}\left| \varphi  \right|^4 + \frac{1} {2}\left\langle {\mu \cdot
\varphi ,\varphi } \right\rangle  - 2\left\langle {F_A^ + ,\mu }
\right\rangle  + \left| \mu  \right|^2},
\end{equation}
where $\mu$ is some fixed imaginary-valued self-dual 2-form and
${\mu  \cdot \varphi }$ represents Clifford multiplication. Then,
we define the perturbed flow equations to be
\begin{equation} \label{flow1pert}
\frac{{\partial \varphi }} {{\partial t}} =  - \nabla _A^* \nabla
_A \varphi  - \frac{1}{4}\left[ {S + \left| \varphi  \right|^2 }
\right]\varphi - \frac{1}{2}\mu  \cdot \varphi,
\end{equation}
\begin{equation} \label{flow2pert}
\frac{{\partial A}} {{\partial t}} =  - d^* F_A   -
i\operatorname{Im} \left\langle {\nabla_A \varphi ,\varphi }
\right\rangle + d^* \mu.
\end{equation}
The purpose of this section is to show that our global existence
and convergence results extend to these perturbed equations.
Rather than duplicate each proof, we will simply outline the
differences. In Lemma \ref{phibound}, we have instead the equation
\[
\frac{\partial } {{\partial t}}\left| \varphi  \right|^2  = -
\Delta \left| \varphi  \right|^2  - 2\left| {\nabla _A \varphi }
\right|^2  - \frac{1} {2}\left[ {S + \left| \varphi \right|^2 }
\right]\left| \varphi  \right|^2 - \operatorname{Re} \left\langle
{\mu  \cdot \varphi ,\varphi } \right\rangle,
\]
where the additional term satisfies $ - \operatorname{Re}
\left\langle {\mu  \cdot \varphi ,\varphi } \right\rangle
\leqslant \left( {\mathop {\max }\limits_{x \in M} \left| {\mu
(x)} \right|} \right)\left| \varphi  \right|^2$. Since
\[
- \frac{1} {2}\left[ {S + {\mathop {\max }\limits_{x \in M} \left|
{\mu (x)} \right|} + \left| \varphi \right|^2 } \right]\left|
\varphi  \right|^2 \leq 0
\]
for $\left| \varphi  \right|$ sufficiently large, the same
argument as before yields a uniform bound on $\sup \left\{ {\left|
{\varphi (x,t)} \right|:x \in M} \right\}$. The other estimates in
section 2 also continue to hold. For the proof of local existence
in section 3, we note that the additional terms are zeroth order
and do not change the parabolicity of the gauge transformed
equations. In section 4, in the proof of Lemma
\ref{firstderivativeestimate}, we have additional terms of
\[
2\operatorname{Re} \left\langle { - \frac{1} {2}\nabla _A (\mu
\cdot \varphi ),\nabla _A \varphi } \right\rangle  \leqslant
c\left| {\nabla _A \varphi } \right|^2  + c\left| {\nabla _A
\varphi } \right|,
\]
and
\[
2\left\langle {dd^* \mu ,dA} \right\rangle  \leqslant c\left| {F_A
} \right|,
\]
and the lemma continues to hold. The proof of Lemma
\ref{regularitytheorem} relies only on Lemma
\ref{firstderivativeestimate}, and is unchanged. For Lemma
\ref{localenergyestimate}, noting that
\[
 - 2\int_M {\phi ^2 \operatorname{Re} \left\langle {\frac{{\partial \varphi }}
{{\partial t}},\frac{1}{2}\mu  \cdot \varphi } \right\rangle }  =
- \frac{d} {{dt}}\int_M {\phi ^2 \frac{1}{2}\left\langle {\mu
\cdot \varphi ,\varphi } \right\rangle },
\]
the proof is entirely analogous. Lemma
\ref{higherderivativeestimate} continues to hold for the same
reason as Lemma \ref{firstderivativeestimate}, as does its
corollary. In Lemma \ref{blowupargument}, the new terms in
(\ref{middleenergy}) are multiplied by factors of $R_m$, and
become negligible in the limit. This establishes global existence.
In Section 5, the proofs of Lemmas \ref{subsequence} and \ref{continuousdependence} are unchanged. In
Lemma \ref{Lojasiewicz}, as for local existence, the additional terms are of
order zero and do not affect parabolicity. Finally, in Lemma \ref{interior},
the additional terms lead to an equation of the same form. The
remaining arguments in this section are unchanged. Thus the analogues of Theorems 1
and 2 hold also for the perturbed equations (\ref{flow1pert}) and
(\ref{flow2pert}), for an arbitrary perturbation parameter $\mu$.

\begin{acknowledgement}
 { The research  of the first author was supported by the Australian Research Council
grant DP0985624.   Theorem 2 was suggested by Prof. Huai-Dong Cao
and Prof. Gang Tian.  We would also like to thank Prof. Huai-Dong
Cao and Prof. Gang Tian for their useful suggestions.}
\end{acknowledgement}

\bibliographystyle{plain}

\end{document}